\documentclass[11pt]{amsart}
\usepackage{amsmath,amscd,latexsym,verbatim,amssymb}
\usepackage[frenchb]{babel}     
\usepackage{times}

\newcommand{\resp}{{\it resp.} }
\newcommand{\cf}{{\it cf.} }
\newcommand{\ie}{{\it i.e.} }
\newcommand{\eg}{{\it e.g.} }

\newcommand{\loccit}{{\it loc. cit.} }
\newcommand{\If}{\Longrightarrow}

\newcommand{\N}{\mathbf{N}}
\newcommand{\Q}{\mathbf{Q}}
\newcommand{\R}{\mathbf{R}}
\newcommand{\boP}{\mathbf{P}}

\newcommand{\Z}{\mathbf{Z}}
\newcommand{\bG}{\mathbb{G}}

\newcommand{\sC}{{\mathcal{C}}}
\newcommand{\sE}{{\mathcal{E}}}

\newcommand{\sI}{{\mathcal{I}}}
\newcommand{\sL}{{\mathcal{L}}}

\newcommand{\sP}{{\mathcal{P}}}

\newcommand{\sR}{{\mathcal{R}}}
\newcommand{\sT}{{\mathcal{T}}}

\newcommand{\inj}{\hookrightarrow}

\newcommand{\surj}{\rightarrow\!\!\!\!\!\rightarrow}

 \newcommand{\un}{{\bf 1}}

\renewcommand{\epsilon}{\varepsilon}
\renewcommand{\phi}{\varphi}

\renewcommand{\lim}{\varprojlim}
\newcommand{\colim}{\varinjlim}

\newcommand{\Ker}{\operatorname{Ker}}

\newcommand{\car}{\operatorname{car}}
\newcounter{spec}
{\end{list}}

\swapnumbers

\newtheorem{thm}{Th\'eor\`eme}[subsection]
\newtheorem{lemme}[thm]{Lemme}
\newtheorem{prop}[thm]{Proposition}
\newtheorem{cor}[thm]{Corollaire}

\newtheorem{sorite}[thm]{Sorite}
\newtheorem{scolie}[thm]{Scolie}
\theoremstyle{definition}
\newtheorem{defn}[thm]{D\'efinition}
\newtheorem{nota}[thm]{Notation}
\newtheorem{ex}[thm]{Exemple}

\newtheorem{rem}[thm]{Remarque}
\newtheorem{rems}[thm]{Remarques}

\newtheorem{compl}[thm]{Compl\'ement}

\numberwithin{equation}{section}

\font\sm=cmr10 at10pt

\newcommand{\prf}{\noindent {\bf D\'emonstration. }}
\renewcommand{\qed}{\hfill $\square$\medskip}

\begin{document}
 \title
{Filtrations de type Hasse-Arf et monodromie $p$-adique}
\author{Yves Andr\'e}
\address{Institut de Math\'ematiques de
Jussieu\\175 rue du
Chevaleret\\75013
Paris\\France.}
\email{andre@math.jussieu.fr }
\rightline{\sm to appear in Invent. Math. (2002)} \maketitle
\tableofcontents

\section*{Introduction} Le principal propos de cet article est d'appliquer la th\'eorie de Christol-Mebkhout \`a
la d\'emonstration de la ``conjecture de monodromie locale $p$-adique" de Crew. 
  
\subsection{} Si $k$ est un corps fini de caract\'eristique $p$, la structure du groupe de Galois $G_{k((z))}=
Gal(k((z))^{sep}/k((z)))$ est plus compliqu\'ee qu'en carac\-t\'eristique $0$, et partant, la th\'eorie de ses
repr\'esentations en est plus riche. 
 
Dans le cas des repr\'esentations $\ell$-adiques (continues, de dimension finie, avec $\ell\neq p$) ou, ce qui
revient au m\^eme, des faisceaux \'etales $\ell$-adiques sur $ Spec\,k((z))$, le th\'eor\`eme de monodromie locale de
Grothendieck affirme que l'inertie $\sI =G_{\bar k((z))}\;$ agit de mani\`ere quasi-unipotente: un
sous\--groupe d'indice fini agit de mani\`ere unipotente; en d'autres termes, la repr\'esentation de l'inertie
devient extension it\'er\'ee de repr\'esentations tri\-viales si on remplace $\bar k((z))$ par une extension finie
s\'eparable convenable (en g\'en\'eral non ab\'elienne). 

\medskip Dans \cite{[Cr2]}, Crew conjecture un analogue $p$-adique de ce th\'eor\`eme, pour les $F$-isocristaux
surconvergents sur $ Spec\,k((z))$, et met en valeur son importance.  

Ces isocristaux surconvergents ``locaux" sont, en langage plus concret (mais qui obs\-curcit un peu l'analogie
avec les faisceaux \'etales $\ell$-adiques), des modules diff\'e\-ren\-tiels $M$ sur le corps diff\'erentiel
$\sE^\dagger$ des fonctions analytiques born\'ees sur une ``mince" couronne $\sC(]r,1[)$ de rayon int\'erieur $r$
non pr\'ecis\'e,
\`a coefficients dans un corps $K$ complet pour une valuation discr\`ete $p$-adique de corps r\'esiduel $k$.   

Un point technique important dans la conjecture de Crew est le passage de $\sE^\dagger$ \`a l'anneau $\sR$ des
fonctions analytiques non n\'ecessairement born\'ees sur une couronne $\sC(]r,1[)$ de rayon int\'erieur non
pr\'ecis\'e (qui a l'avantage sur $\sE^\dagger$ d'\^etre auto-dual pour la dualit\'e locale de Dwork). 

Toute extension finie s\'eparable de $k((z))$ donne naissance \`a
une extension finie de
$\sE^\dagger$ et de $\sR$ respectivement, et la conjecture affirme que si $M$ admet une structure de Frobenius,
alors $M$ est {\it quasi-unipotent}, \ie devient extension it\'er\'ee de modules diff\'erentiels triviaux sur une
extension finie de $\sR$ provenant d'une extension finie s\'eparable de $k((z))$. Cette conjecture ``explique"
les isocristaux surconvergents sur $ Spec\,k((z))$ admettant une structure de Frobenius en termes des
repr\'esentations $p$-adiques de
$G_{k((z))}$ dont l'image de l'inertie est finie.   

\medskip 
Le passage de $\sE^\dagger$ (en faveur chez Dwork-Robba) \`a $\sR$ \'etait aussi la cl\'e des progr\`es
que Christol et Mebkhout ont accomplis dans l'\'etude des \'equa\-tions diff\'eren\-tiel\-les $p$-adiques,
avec pour motivation la conjecture de l'indice de Robba. L'une de leurs innovations majeures est la th\'eorie des
``pentes $p$-adiques" des modules diff\'erentiels sur $\sR$, d\'efinies en terme de convergence de solutions en
divers points g\'en\'eriques. Comme application de leur th\'eor\`eme de l'indice, ils obtiennent une propri\'et\'e
d'int\'egralit\'e \`a la Hasse-Arf pour la filtration par les pentes $p$-adiques. 

 C'est l\`a l'outil essentiel dont nous nous servons pour d\'emontrer la conjecture de Crew, sous la forme
renforc\'ee suivante (et sans supposer $k$ fini).

\begin{thm}\label{mono} Tout module diff\'erentiel de pr\'esentation finie sur $\sR$ admettant une structure de
Frobenius est quasi-unipotent, \ie poss\`ede une base de solutions dans $\sR'[\log z]$, o\`u $\sR'$ est l'extension
finie
\'etale de
$\sR$ issue d'une extension finie s\'eparable convenable de $k((z))$. 
\end{thm} 

Notre r\'esultat est plus pr\'ecis: il d\'etermine le groupe tannakien attach\'e \`a la cat\'egorie des
modules diff\'erentiels admettant une structure de
Frobenius sur une ``couronne $p$-adique infiniment mince" (th\'eor\`eme \ref{main}), et relie directement cette
cat\'egorie \`a celle des repr\'esentations $p$-adiques du corps local $k((z))$ \`a inertie finie.

En fait, nous \'elucidons d'abord la structure de cette cat\'egorie de modules dif\-f\'e\-ren\-tiels, pour en
d\'eduire la structure de ses objets.

\subsection{}  On rencontre deux types
g\'en\'eraux de filtrations dans diverses situations tannakiennes concr\`etes: celles du type ``pentes
frobeniusiennes", qui se comportent ``additivement" pour le produit tensoriel, et qui ont \'et\'e formalis\'ees par
Saavedra \cite[IV]{saavedra}; et celles pour lesquelles le produit tensoriel de deux objets de
pentes $\leq \lambda\,$ est de pentes $\,\leq \lambda$ (c'est le cas de la filtration par les pentes $p$-adiques). 

Nous formalisons ce second type de filtration par les pentes, dans une cat\'egorie
tannakienne quelconque. On peut attacher \`a tout objet son poly\-g\^one de Newton suivant la recette
classique. Dans tous les exemples int\'eres\-sants, il s'av\`ere que le poly\-g\^one de Newton de tout objet  
est \`a sommets entiers. Nous appelons {\it filtration de type Hasse-Arf} une filtration par les pentes ayant
cette propri\'et\'e remarquable.

 En pratique, comme l'illustre l'exemple originel de la filtration de Hasse-Arf sur les
repr\'esentations finies d'un groupe de Galois local, on dispose en outre de
fonctorialit\'es (changements de base finis) v\'erifiant certaines compatibilit\'es partielles (de type
Herbrand) vis-\`a-vis des filtrations, qui enrichissent beaucoup la th\'eorie.  Nous adaptons une part de ce
formalisme au cas des filtrations de type Hasse-Arf abstraites. 

On obtient alors une situation o\`u interf\`erent la th\'eorie des
repr\'esentations (aspect tannakien), la combinatoire (int\'egralit\'e des polyg\^o\-nes de Newton), et
des transferts. Cette situation est extr\^emement rigide. On peut d\'emontrer, sous certaines
hypoth\`eses sur les objets de dimension un, un th\'eor\`eme
 de structure g\'en\'eral (\ref{struc}) qui est la principale innovation technique de l'article.  
   
Ce th\'eor\`eme n'a rien de ``$p$-adique": la possibilit\'e de l'appliquer aux isocristaux
surconvergents repose essentiellement sur les r\'esultats de Christol-Mebkhout. 
{\it A posteriori}, le th\'eor\`eme de monodromie locale $p$-adique \ref{mono} permet de retrouver et
comprendre beaucoup de leurs \'enonc\'es de mani\`ere galoisienne, et m\^eme d'aller
plus loin dans l'ana\-lyse des modules diff\'eren\-tiels irr\'eductibles (\S\ref {s8}).

\subsection{} Le th\'eor\`eme de monodromie locale $\ell$-adique de Grothendieck vaut non seulement pour les
corps locaux de caract\'eristique $p>0$, mais aussi pour les corps locaux $K$ d'in\'egales caract\'eristiques
$(0,p)$. Pour un tel corps, l'analogue $p$-adique du th\'eor\`eme de monodromie locale est la {\it conjecture de
Fontaine} \cite{[F1]} sur les repr\'esentations $p$-adiques de corps $p$-adiques, qui s'\'enonce ainsi: 

\medskip {\it Toute repr\'esentation de De Rham est potentiellement semistable.}  

\medskip\noindent Le lien entre cette conjecture et celle de Crew a \'et\'e mis en lumi\`ere dans un travail
r\'ecent remarquable de Berger \cite{ber} (un pas ant\'erieur important \'etant le th\'eor\`eme de surconvergence de
\cite{[ChCo]}). Il y d\'emontre entre autres: 

\begin{thm}\label{berger} $ \rm ( Berger \,$\cite[5]{ber} $\rm )$. Il existe un $\otimes$-foncteur fid\`ele et exact
de la cat\'egorie des repr\'esentations $p$-adiques de De Rham de $K$ vers celle des modules diff\'erentiels de
pr\'esentation finie sur
$\sR_K$ admettant une structure de Frobenius. En outre, une repr\'esentation est potentiellement semistable si et
seulement si le module diff\'erentiel associ\'e est quasi-unipotent.  
\end{thm}  

La conjecture de Fontaine r\'esulte imm\'ediatement des th\'eor\`emes \ref{berger} et \ref{mono}. 

\subsection{} L'article se termine par une br\`eve discussion de la conjecture de Dwork sur les structures de
Frobenius, et des extensions canoniques.
 
\bigskip\bigskip {\sm Les id\'ees contenues dans ce texte sont une \'evolution de celles de \cite{an} (o\`u un
contre-exemple putatif propos\'e par Z. Mebkhout, longtemps rest\'e en suspens, \'etait mis au ban d'essai). Je
remercie Bruno Chiarellotto et Pierre Colmez de m'avoir fortement encourag\'e
\`a les r\'ediger enfin, en m'expliquant le lien avec le travail de L. Berger et la conjecture de J.-M. Fontaine. Je
les remercie aussi, ainsi que R. Crew, pour leurs commentaires sur le texte.}
 \newpage

\section{Filtrations par les pentes dans les cat\'egories tannakiennes}\label{s1} 

\subsection{} On formalise ici quelques propri\'et\'es des filtrations par les pentes rencontr\'ees dans divers
contextes (voir \cite[II]{[K1]} pour le cadre diff\'erentiel, \cite[I]{[K2]} pour le cadre
des re\-pr\'e\-sen\-ta\-tions galoisiennes $\ell$-adiques). 

\medskip On se donne une cat\'egorie tannakienne $\sT$ sur un corps ${E}$. Par {\it sous-cat\'egorie tannakienne} de
$\sT$, on entendra toujours une sous-cat\'egorie strictement pleine contenant l'unit\'e $\un$, stable par
sous-quotient (au sens de $\sT$), par $\oplus$, $\otimes$ et dualit\'e ${}^\vee$.

 On se donne un syst\`eme d\'ecroissant, index\'e par $\lambda \in \R_{\geq 0}$, de {\it sous-foncteurs
${E}$-lin\'eaires exacts}
$F^{>\lambda}$ du foncteur identique de
$\sT$.   
 
Ces endofoncteurs d\'efinissent un syst\`eme d'endofoncteurs ${E}$-lin\'eaires exacts \[\displaystyle{gr^\lambda
=(\lim_{\mu< \lambda}\, F^{>\mu})/F^{>\lambda}}\;\hbox{ pour}\;  \lambda>0 , \;\hbox{et}  \;\displaystyle{gr^0
=id/F^{>0}}\;\;\] (pour v\'erifier l'exactitude, on peut
remarquer que comme tout objet
$M$ est de longueur finie, on a 
 $\lim_{\mu< \lambda}\, F^{>\mu}(M)= F^{>\nu}(M)$ pour $\nu<\lambda$ conve\-nable).

 On note $\sT_{\leq \lambda}$ la sous-cat\'egorie pleine de $\sT$ form\'ee des objets sur lesquels $ 
F^{>\lambda}$ s'annule. C'est une cat\'egorie ab\'elienne, en vertu de l'exactitude de $F^{>\lambda}$.

Faisons en outre les hypoth\`eses suivantes: pour tous $M,N\in Ob(\sT)$ et tout $\lambda\geq 0$, on a 

\medskip
$FP_1)$   $\;\;\displaystyle\lim_{\lambda \geq 0} \; F^{>\lambda}$ est le foncteur nul,

$FP_2)$    $\;\;\displaystyle\colim_{\mu> \lambda} \; F^{>\mu}= F^{>\lambda}$,

$FP_3)$   $\;\;F^{>\lambda}(\un)=0$, 

$FP_4)$   $\;\;F^{>\lambda}(M)= F^{>\lambda}(N)= 0 \If F^{>\lambda}(M\otimes N)=0.$

$FP_5)$   $\;\;F^{>\lambda}(M)=  0 \If  F^{>\lambda}(M^\vee)=0.$

\medskip
\noindent Alors les $\sT_{\leq \lambda}$ forment un syst\`eme croissant de sous-cat\'egories
tannakiennes de $\sT$, dont la r\'eunion est \'egale \`a $\sT$. Notons que la condition de ``continuit\'e \`a droite"
$FP_2)$ se traduit par:
\[\displaystyle{\sT_{\leq \lambda}= \;\bigcap_{\mu >\lambda }  \;\sT_{\leq\mu}.}\]

\begin{defn} Un syst\`eme de foncteurs $F^{>\lambda}$ comme ci-dessus sera appel\'e une \emph{filtration par les
pentes} de $\sT$ \footnote{on prendra garde \`a ne pas confondre cette notion avec celle de $\otimes$-filtration
introduite dans \cite[IV.2.1]{saavedra}, qui n'a rien \`a voir.}.
\\ Comme tout objet $M$ de $\sT$ est de longueur finie, sa filtration
$F^{>\lambda}(M)$ n'a qu'un nombre fini de sauts $\lambda_0,\lambda_1,\ldots $  (tels que $gr_{\lambda_i}(M) \neq
0)$, qui d\'efinissent les {\it pentes} de $M$. Il est parfois utile de convenir que $0$ est de pente $-\infty$.
\end{defn}

\begin{rem} Les conditions $FP_4)$  et
$FP_5)$ se traduisent par: si les pentes de $M$ et $N$ sont $\leq \lambda$, il en est de m\^eme des pentes de
$M\otimes N,M^\vee$ et $N^\vee$.
\end{rem}

\begin{lemme}\label{dec} On a une d\'ecomposition canonique, fonctorielle, de tout objet $M=\oplus gr^\mu(M)$,
et $F^{>\lambda}(M)=\oplus_{\mu>\lambda}\,gr^\mu(M)$. Pour $\mu\neq \lambda$, on a
$Hom_\sT(gr^\lambda(M),gr^\mu(M))=0$.
\end{lemme}

\prf De la condition $FP_5)$, on tire, par sym\'etrie, que $gr^\lambda(M)=0\iff
gr^\lambda(M^\vee)=0$. Soit $\lambda_0$ la plus petite pente de $M$. Par r\'ecurrence sur la longueur de $M$,
il suffit, pour la premi\`ere assertion, de montrer que la projection $M\surj gr^{\lambda_0}(M)$ admet une section. 

Consid\'erons le monomorphisme $gr^{\lambda_0}(M)^\vee\inj M^\vee$. Son image est purement de pente $\lambda_0$,
donc le compos\'e $gr^{\lambda_0}(M)^\vee\inj M^\vee\surj gr^{\lambda_0}(M^\vee)$ est encore un monomorphisme.
En passant au dual, on trouve que le compos\'e $(gr^{\lambda_0}(M^\vee))^\vee \inj M \surj gr^{\lambda_0}(M)$
est un \'epimorphisme, d'o\`u la premi\`ere assertion. 

La seconde assertion r\'esulte de l\`a et de ce que $gr^\lambda gr^\mu(M)=0$. 
 \qed

\begin{cor}\label{cor ind} Tout objet ind\'ecomposable de $\sT$ n'a qu'une seule pente.\qed 
\end{cor}

\subsection{} On suppose maintenant $\sT$ munie d'un foncteur fibre 
\[\omega: \sT\to Vec_{E}.\]
Notons $G$ le ${E}$-groupe affine $Aut^\otimes \omega$. Par la th\'eorie tannakienne, le syst\`eme des $\sT_{\leq
\lambda}$ d\'efinit un syst\`eme d\'ecroissant de sous-groupe ferm\'es normaux $G^{>\lambda}\triangleleft G$, tels
que $G_{\leq \lambda}:= G/G^{>\lambda} \cong  Aut^\otimes \omega_{\vert \sT_{\leq \lambda}}.$ Pour
$\lambda = 0 $, on \'ecrit plut\^ot $G_0$ que $G_{\leq 0}$ (puisqu'il n'y a pas de pente n\'egative). Terminologie:
$G^{>0}$ est le ``{\it sous-groupe sauvage}", et $G_0$ le ``{\it quotient mod\'er\'e}".

 On a
  \[\omega(M)^{G^{>\lambda}}=\oplus_{\mu \leq \lambda} \;\omega(gr^\mu(M))\] 
\[ \bigcap
G^{>\lambda}= 1 \;\;(\hbox{par}\; FP_1) ,\] 
tandis que la condition $FP_2)$ se traduit par:
\\$SC)$ $\;\;\;\;\displaystyle{G^{>\lambda}= \;\overline{\bigcup_{\mu >\lambda }  \;G^{>\mu}}\;\;\;\;
\hbox{(adh\'erence de Zariski de la r\'eunion)}.}$

\begin{rems}\label{sous/ext} a) Si $\sT'\subset \sT$ est une sous-cat\'egorie tannakienne, la restriction
\`a
$\sT'$ de la filtration par les pentes de $\sT$ d\'efinit une filtration de $Aut^\otimes \omega_{\vert \sT'}$ qui
n'est autre que la filtration par l'image des $G^{>\lambda}$ par l'\'epi\-mor\-phis\-me canonique $G=Aut^\otimes
\omega\to Aut^\otimes
\omega_{\vert \sT'}\;\;$  ($\sT'_{\leq \lambda}= \sT_{\leq \lambda}\cap \sT'$). 

b) Soit ${E}'/{E}$ une extension finie. On a une notion d'extension des scalaires $\sT_{({E}')}$ \`a la Saavedra \cf
\cite[III.1]{saavedra}: c'est la cat\'egorie des ${E}'$-modules dans $\sT\;$\footnote{si ${E}'/{E}$ est s\'eparable,
$\sT_{({E}')}$ peut aussi se d\'ecrire comme l'enveloppe pseudo-ab\'elienne de la cat\'egorie obtenue \`a partir de
$\sT$ en \'etendant les scalaires \`a ${E}'\;$, \cf\cite[5.3.2]{ak}.};
$\sT_{({E}')}$ est naturellement
\'equivalente \`a $Rep_{{E}'}\,(G\otimes_{E} {{E}'})$, et on a un foncteur $M\mapsto M\otimes {E}'$ de $\sT$ vers
$\sT_{({E}')}$. Si $\sT$ est munie d'une filtration par les pentes, il existe une unique filtration par les pentes
$\sT_{({E}')}$ compatible avec ce foncteur: les pentes d'un ${E}'$-module dans $\sT$ sont les pentes de l'objet de
$\sT$ sous-jacent. Les groupes $G_{{E}'}^{>\lambda}$ correspondants ne sont autres que les $G^{>\lambda}\otimes_{E}
{E}'$.    
\end{rems}
 
 \begin{nota}\label{not} Pour tout $\lambda >0$, on pose $\;\;\displaystyle{G^{(\lambda)}=
\bigcap_{\mu<\lambda} \, G^{>\mu} }.$
  \end{nota}
\noindent Ce sont des sous-groupes ferm\'es normaux de $G$, et on a les propri\'et\'es sui\-vantes: 
\[\displaystyle{\bigcap_{\lambda >0}\, G^{(\lambda)}= 1,\;\;\; G^{>\lambda}=\overline{\bigcup_{\mu >\lambda}\;
G^{(\mu)}} \;\;\;\forall \lambda \geq 0 ,}\;\; \; G^{>\lambda}\subset G^{(\lambda)} \;\,\forall \lambda > 0, \]
\\$CG)$ $\displaystyle{\;\;\;\;\;\bigcap_{\lambda >\mu>0}\, G^{(\mu)}=
G^{(\lambda)}\;\;}$  (continuit\'e \`a gauche), 

\medskip\noindent $SS1)\;$ pour tout $M\in Ob(\sT)\;$ et tout $\lambda>0$, la sous-re\-pr\'e\-sen\-ta\-tion triviale
$\omega(M)^{G^{(\lambda)}}$ est facteur direct de $\omega(M)$ en tant que re\-pr\'e\-sen\-ta\-tion de $G^{(\lambda)}$.

\medskip Cette derni\`ere condition (``semisimplicit\'e pour la valeur propre $1$") vient de ce que pour tout
$\lambda >0$, \[\omega(M)^{G^{(\lambda)}}=\oplus_{\mu < \lambda} \;\omega(gr^\mu(M))\] admet 
  $\; \oplus_{\mu \geq \lambda}\;\omega( gr^\mu(M))$ comme suppl\'ementaire stable sous $G^{(\lambda)}$.

\begin{thm}\label{P2} La donn\'ee d'une filtration par les pentes $\;(F^{>\lambda})_{\lambda\geq 0}\;$ sur $\sT$
\'equivaut
\`a celle d'une filtration d\'ecroissante s\'epar\'ee $\;(G^{(\lambda)}\triangleleft G)_{\lambda >0}$ par des
sous-groupes ferm\'es normaux v\'erifiant les conditions $CG)$ et $SS1)$.
\end{thm}   

\prf Pour toute filtration par les pentes, on vient de voir que syst\`eme des $G^{(\lambda)}$ d\'efinis en
\ref{not} v\'erifient les conditions du th\'eor\`eme. Observons en outre que $F^{>\lambda}(M)$ est
d\'etermin\'e par 
\[\omega(F^{>\lambda}(M))= \Ker(\omega(M)\to \omega(M)_{G^{>\lambda}})\]
o\`u $\omega(M)_{G^{(\lambda)}}$ d\'esigne l'espace quotient des co-invariants (qui est cano\-ni\-quement isomorphe
\`a l'espace des invariants puisque ce der\-nier admet un suppl\'ementaire stable). 

Passons \`a la r\'eciproque. On se donne une filtration $\;G^{(\lambda)}\triangleleft G\;$ comme dans le
th\'eor\`eme. Pour tout
$\lambda\geq 0$, on pose $G^{>\lambda}=\overline{\bigcup_{\mu >\lambda}\; G^{(\mu)}}$. Ils forment aussi une
filtration d\'ecrois\-sante s\'epar\'ee par des sous-groupes ferm\'es normaux, et v\'erifient l'analogue de la
condition $SS1)$ (en fait l'image de $G^{>\lambda}$ dans $GL(\omega(M))$ s'identifie
\`a l'image de l'un des $G^{(\mu)}$). 
Puisque $G^{>\lambda}$ est normal dans $G$, $ \Ker(\omega(M)\to \omega(M)_{G^{>\lambda}})$ est $G$-stable, donc 
s'identifie \`a $\omega(F^{>\lambda}(M))$ pour un sous-objet
$F^{>\lambda}(M)$ de $M$. On obtient ainsi un foncteur $M\mapsto
F^{>\lambda}(M)$. 

Compte tenu de $SS1)$, la formation des co-invariants est exacte, donc $M\mapsto
\omega(F^{>\lambda}(M))$ est exact, et finalement  $M\mapsto F^{>\lambda}(M)$ est exact puisque $\omega$ est
fid\`ele et exact. Les conditions $FP_1)$ et $FP_3)$ sont clairement remplies. La condition $FP_2)$ se traduit par
$SC)$, qui est avec la d\'efinition ci-dessus de $G^{>\lambda}$ une tautologie: $\overline{\bigcup_{\nu
>\lambda}\; G^{(\nu)}}=\overline{\bigcup_{\mu
>\lambda} \overline{\bigcup_{\nu
>\mu}\,G^{(\nu)}}}$. Quant \`a $FP_4)$ et $FP_5)$, elles suivent de ce que la sous-cat\'egorie pleine de
$\sT$ form\'ee des objets $M$ tels que $G^{>\lambda}$ agisse trivialement sur $\omega(M)$ est une sous-cat\'egorie
tannakienne, ce qui est clair.

Il est aussi clair que cette construction $(G^{(\lambda)})\mapsto (F^{>\lambda})$ est inverse \`a gauche de la
construction $ (F^{>\lambda})\mapsto (G^{(\lambda)}) $ de \ref{not}. Il reste \`a voir que c'est aussi un inverse
\`a droite, ce qui revient \`a voir que 
\[\displaystyle{G^{(\lambda)}= \bigcap_{\mu<\lambda} \overline{\bigcup_{\nu >\mu}  \; G^{(\nu)}}}. \]
L'inclusion $\subset$ est claire; il s'agit de montrer l'inclusion oppos\'ee.  
C'est l\`a qu'intervient la condition $CG)$: $\;\displaystyle{G^{(\lambda)}=
\bigcap_{\mu<\lambda}\, G^{(\mu)}}$. On conclut du fait que $ G^{(\mu)}\supset \overline{\bigcup_{\nu >\mu}  \;
G^{(\nu)}}$.
\qed
 
\begin{cor} On a les inclusions:
  \[\displaystyle{\oplus_{\mu<\lambda}\;gr^\lambda(M)\otimes gr^\mu(N)\subset gr^\lambda(M\otimes N),
}\] 
\[\displaystyle{gr^\lambda(M)\otimes gr^\lambda(N)\subset \oplus_{\mu \leq\lambda}\;gr^\mu(M\otimes N), 
}\]
 \[\displaystyle{\oplus_{\mu<\lambda}\;(\underline{Hom}(gr^\lambda(M), gr^\mu(N)\oplus \underline{Hom}(gr^\mu(M),
gr^\lambda(N)))\subset gr^\lambda(\underline{Hom}(M , N)),} \] 
 \[\displaystyle{\underline{Hom}(gr^\lambda(M), gr^\lambda(N))\subset \oplus_{\mu \leq
\lambda}\;gr^\mu(\underline{Hom}(M , N)). } \] 
\end{cor}
 
\prf Par $FP_4)$ et $FP_5)$, on a \[\displaystyle{gr^\lambda(M)\otimes gr^\mu(N)\subset \oplus_{\nu \leq
max(\lambda,\mu)}\;gr^\nu(M\otimes N) ,  }\]
\[\displaystyle{\underline{Hom}(gr^\lambda(M), gr^\mu(N))\subset \oplus_{\nu \leq
max(\lambda,\mu)}\;gr^\nu(\underline{Hom}(M , N)), } \] d'o\`u les deuxi\`eme et quatri\`eme inclusions du
corollaire. D\'emontrons la premi\`ere (la troisi\`eme est analogue). On a
$\;(\omega(gr^\lambda(M)))^{G^{(\lambda)}}=0\;$ et
$\;(\omega(gr^\mu(N)))^{G^{(\lambda)}}=\omega(gr^\mu(N))\;$ si
$\mu<\lambda$. On a donc $\;(\omega(gr^\lambda(M)\otimes gr^\mu(N)))^{G^{(\lambda)}}=0$, d'o\`u
$\; gr^\lambda(M)\otimes gr^\mu(N)= gr^\lambda(  gr^\lambda(M)\otimes gr^\mu(N))$.\qed
 
\begin{rem}\label{NB} Les $G^{(\lambda)}$ \'etant des sous-groupes normaux non n\'eces\-sai\-rement
caract\'eristiques, on doit prendre garde au fait que la filtration par les pentes n'est pas n\'ecessairement
invariante par auto-\'equivalence de $\sT$. On a toutefois le

\begin{sorite}\label{sor} Toute auto-\'equivalence de $\sT$ qui est naturellement isomorphe au $\otimes$-foncteur
identique respecte la filtration par les pentes. En particulier, si $u:\sT\to \sT'$ est une \'equivalence de
cat\'egories tannakiennes, et si on munit $\sT'$ de la filtration par les pentes induite par $u$, tout quasi-inverse
respecte les filtrations par les pentes.    
\end{sorite}
\end{rem}

\section{Polyg\^ones de Newton et filtrations de Hasse-Arf}\label{s2} 

\subsection{Polyg\^ones de Newton} Soit $\sT$ une cat\'egorie tannakienne sur ${E}$ munie d'une filtration par les
pentes
$(F^{>\lambda})$. Soit $M$ un objet de $\sT$. Soient comme ci-dessus
$\lambda_i, i=0,1\ldots,$ les pentes de
$M$, rang\'ees dans l'ordre croissant. 

On d\'efinit le {\it polyg\^one de Newton} $NP(M)$ comme l'enveloppe convexe dans
$\R_{\geq 0}^2$ des points $(\sum_0^i \dim gr^{\lambda_j}(M), \sum_0^i {\lambda_j}\dim gr^{\lambda_j}(M)), i\geq 0$,
auxquels on joint par convention l'axe vertical. 

L'ordonn\'ee du sommet le plus \`a droite (\ie d'abscisse maximale) de $NP(M)$ s'appelle la {\it hauteur} du
polyg\^one de Newton.
 
\begin{lemme}\label{Npol} La r\`egle $M\mapsto NP(M)$ s'\'etend en un {\it homomorphisme } $NP$ du groupe de
Grothen\-dieck
$K_0(\sT)$ (pour les suites exactes) vers le groupe attach\'e au mono\"{\i}de additif des sous-ensembles convexes
poly\-g\^o\-naux de
$\R_{\geq 0}^2$.  En outre $NP(M^\vee)=NP(M)$.    
\end{lemme} 

C'est clair, compte tenu de l'exactitude de $gr^\lambda$ et de la propri\'et\'e $FP_5)$ des filtrations par
les pentes. \qed

\begin{lemme}\label{hauteur} La r\`egle $M\mapsto hauteur(NP(M))$ s'\'etend en un homomorphisme
\[ hNP : K_0(\sT)\to \R\] qui d\'etermine compl\`etement la filtration par les pentes. 
\end{lemme} 

 Noter que par \ref{cor ind}, tout ind\'ecomposable $M$ a une seule pente, donn\'ee par $hNP(M)/\dim(M)$. Par
Krull-Schmidt, il est alors clair que $hNP$ d\'etermine compl\`etement la filtration par les pentes.\qed

\subsection{Filtrations de Hasse-Arf} \begin{defn} Une {\it filtration de type Hasse-Arf} - ou plus bri\`evement:
filtration de Hasse-Arf - est une filtration par les pentes telle que tous les polyg\^ones de Newton soient {\it \`a
sommets dans
$\N^2$}.
 \end{defn}

Cette propri\'et\'e d'{\it int\'egralit\'e} joue un r\^ole essentiel dans la suite. 
On notera que dans le cas d'une
filtration de Hasse-Arf, les pentes de tout objet sont des nombres rationnels. 

\begin{lemme}\label{int} Supposons donn\'ee une filtration par les pentes $(F^{>\lambda})$ sur $\sT$. Les conditions
suivantes sont
\'equivalentes:
\\ a) $(F^{>\lambda})$ est une filtration de Hasse-Arf,
\\ b) pour tout objet $M$, la hauteur du polyg\^one de Newton de $M$ est un entier. 
\\ c) $hNP$ est \`a valeurs dans $\Z$. 
\\ d) pour tout objet irr\'eductible $M$, le produit de la pente de $M$ par la dimension de $M$ est un entier. 
\end{lemme}

\prf En vertu de \ref{Npol}, $a)$  
\'equivaut \`a l'int\'egralit\'e du polyg\^one de Newton des objets irr\'eductibles, ce qui se traduit
indiff\'eremment par
$b),\,c)$ ou $d)$.
\qed
 
Nous sugg\'erons au lecteur de lire l'appendice A pour voir une manifestation de la rigidit\'e qu'impose cette
condition ``combinatoire"
dans le cas simple o\`u les pentes $\neq 0$ sont non enti\`eres.

 \section{Exemples de filtrations de Hasse-Arf}\label{s3}

\subsection{Re\-pr\'e\-sen\-ta\-tions galoisiennes finies}\label{gal}

Il s'agit de l'exemple de base qui justifie la terminologie pr\'ec\'edente. Soit $K$ est un corps henselien pour une
valuation discr\`ete
\`a corps r\'esiduel
$k$ parfait d'exposant caract\'eristique $p\geq 1$. On a des suites exactes
\[1\to \sI \to G_K= Gal(K^{sep}/K) \to Gal(k^{sep}/k)\to 1 \]
\[1\to \sP \to \sI\to \Pi_{\ell\neq \car k}\,\Z_\ell\to 1 \] 
o\`u $\sP$ est un pro-$p$-groupe (donc trivial si $p=1$). On a $\sI=G_{K^{hs}}$, o\`u $K^{hs}$ est l'henselis\'e
strict de $K$.  

On consid\`ere ici $G_K, \sI$ et $\sP$ comme sch\'emas en groupes profinis constants sur ${E}$, de sorte que
les objets de $Rep_{E}(G_K)$ sont des re\-pr\'e\-sen\-ta\-tions d'ima\-ge finie par d\'efinition.

\medskip
Soient $\;G_K^{(\lambda)} \triangleleft G_K$ les groupes de ramification de $K$ en
num\'erotation sup\'erieure (\cf \cite{[Se]}). Alors par le th\'eor\`eme de Hasse-Arf (et via la proposition
\ref{P2}), les groupes $G_K^{(\lambda)}$  munissent
  $Rep_{E}(G_K)$ d'une filtration
de Hasse-Arf. La hauteur du polyg\^one de Newton d'une re\-pr\'e\-sen\-ta\-tion $V$ de $G_K$ s'appelle le {\it
conducteur de Swan de } $V$, \cf \cite{[K2]}. Le th\'eor\`eme de Hasse-Arf \'equivaut \`a dire que c'est
toujours un entier\footnote{r\'ecemment, Abbes et Saito ont d\'efini une filtration par les pentes sur
$Rep_{E}(G_K)$ sans l'hypoth\`ese que $k$ soit parfait \cite{as}, mais il n'est pas clair que ce soit une
filtration de Hasse-Arf.}.  

Les groupes
$G_K^{>\lambda}$ sont aussi parfois not\'es $G_K^{(\lambda+)}$.  Notons que le groupe d'inertie $G_K^{(0)}$ n'est pas
pris en compte par le formalisme des filtrations par les pentes: les pentes ne distinguent pas entre extensions non
ramifi\'ees et extensions mod\'er\'ement ramifi\'ees.  

\subsection{Re\-pr\'e\-sen\-ta\-tions $\ell$-adiques}\label{ell}

Ici on prend ${E}=\bar\Q_\ell$ ($\ell\neq  p$), et $K$ comme ci-dessus. On consid\`ere la cat\'egorie
tannakienne $Rep_{cont}(G_K/\bar\Q_\ell)$ des re\-pr\'e\-sen\-ta\-tions
$\ell$-adiques du groupe compact $G_K$, c'est-\`a-dire des re\-pr\'e\-sen\-ta\-tions dans un $\bar\Q_\ell$-espace de
dimension finie qui proviennent de re\-pr\'e\-sen\-ta\-tions continues d\'efinies sur une extension finie non
pr\'ecis\'ee de $\Q_\ell$. On note $Rep_{cont}^F(\sI/\bar\Q_\ell)$ la cat\'egorie tannakienne qu'on obtient en
appliquant le $\otimes$-foncteur $Res_{G_K}^\sI$. C'est une sous-cat\'egorie tannakienne de
$Rep_{cont}(\sI/\bar\Q_\ell)$  (si $k$ est fini, on peut la voir comme la sous-cat\'egorie des repr\'esentations 
admettant ``une structure de Frobenius").  

Les groupes de ramification de $G_K$ d\'efinissent une filtration de Hasse-Arf sur $Rep_{cont}^F(\sI/\bar\Q_\ell)$,
\'etudi\'ee en d\'etail dans \cite{[K2]}.

Un argument de compacit\'e bien connu montre que $\sP$ agit toujours \`a travers un groupe fini (\cf
\cite[1.10]{[K2]}). Une version du th\'eor\`eme de mono\-dromie locale $\ell$-adique de Grothendieck
\cite[App.]{[ST]} s'\'enonce ainsi: 
 supposons que 

$(\ast)\;$ aucune
extension finie de
$k$ ne contient toutes les racines de l'unit\'e d'ordre une puissance de
$\ell$. (C'est \'evidemment le cas si $k$  est fini). 
 
\medskip {\it  Alors tout objet de $Rep_{cont}^F(\sI/\bar\Q_\ell)$ est quasi-unipotent}, \ie la restriction de la
repr\'esentation \`a un sous-groupe ouvert de $\sI$ est unipotente (et agit \`a travers le facteur $\Z_\ell$ de
l'inertie mod\'er\'ee).  

\medskip
Si ${}^FG$ est le groupe tannakien sur $\bar\Q_\ell$ associ\'e au foncteur ``espace sous-jacent", de sorte que
$Rep_{cont}^F(\sI/\bar\Q_\ell)\sim Rep_{\bar\Q_\ell}\,{}^FG\;$,  on d\'eduit de ce th\'eor\`eme que  
\[ {}^FG \cong \sI \times \bG_a,\;\;
 {}^FG^{>0}\cong \sP \] 
comme $\bar\Q_\ell$-groupes affines. La filtration $({}^FG^{(\lambda)})$ correspond \'evidemment \`a celle
de $\sI$ par l'isomorphisme pr\'ec\'edent. 

\subsection{Modules diff\'erentiels sur un corps local d'\'egale caract\'eristique nulle}\label{levelt} Soit
${E}((z))$ un corps de s\'eries formelles sur un corps ${E}$ de carac\-t\'eristique nulle. Soit
$MC({E}((z))/{E})$ la cat\'egorie tannakienne sur ${E}$ des ${E}((z))[{d\over dz}]$-modules de ${E}((z))$-dimension
finie. 

La th\'eorie formelle des modules diff\'erentiels (Turrittin, Levelt...) associe fonctoriellement \`a
tout objet de $MC({E}((z))/{E})$ une filtration selon les pentes, \cf \cite[II]{[K1]}. La hauteur du polyg\^one de
Newton d'un objet $M$ de $MC({E}((z))/{E})$ s'appelle l'{\it irr\'egularit\'e (formelle) de } $M$.
C'est toujours un entier (trivialement), qui s'interpr\`ete comme un indice, \cf \cite{[Ma]}. 

On obtient ainsi une filtration de Hasse-Arf sur $MC({E}((z))/{E})$. La structure des groupes $G^{(\lambda)}$ est
\'etudi\'ee en d\'etail dans \cite[II]{[K1]}.

Dans \loccit, Katz construit aussi une extension canonique de tout objet de $MC({E}((z))/{E})$ en un module
diff\'erentiel sur $\boP^1\setminus\{0,\infty\}$, d'o\`u un foncteur fibre canonique $\omega: \,MC({E}((z))/{E})\to
Vec_{E}$ en prenant la fibre au point $1$.

\begin{sloppypar}
\subsection{Modules diff\'erentiels sur des couronnes $p$-adiques}\label{CM} Soit $K$ un corps complet pour une
valuation discr\`ete - ou plus g\'en\'eralement un corps {\it maximalement complet}\footnote{\ie tel que toute
intersection de disques embo\^{\i}t\'es non vides est non vide.} - d'in\'egales ca\-rac\-t\'eristiques. On note
$k$ le corps r\'esiduel (de ca\-rac\-t\'eristique $p>0 $). 

Soit $\sR=\sR_{K,z}$ l'anneau (int\`egre) des fonctions analytiques dans un couronne non circonf\'erenci\'ee
$\sC(]r,1[)$ de diam\`etre int\'erieur non pr\'ecis\'e\footnote{Cet anneau diff\'erentiel appara\^{\i}t dans les
travaux de Robba (\eg \cite[2]{R1}, \cite[3.11]{R2}).}. Soit
$MC(\sR/K)$ la cat\'egorie tannakienne sur
$K$ des $\sR[{d\over dz}]$-modules de pr\'esentation finie sur $\sR$. On peut montrer que tout tel module est en fait
automatiquement {\it libre} de type fini sur $\sR$ (\cite[2]{an}).    

On consid\`ere la sous-cat\'egorie pleine $MCS(\sR/K)$ de $MC(\sR/K)$ for\-m\'ee des modules diff\'erentiels {\it
solubles} vers le bord ext\'erieur, \ie dont le rayon de convergence au point g\'en\'erique de module $\rho$ tend
vers $1$ avec $\rho$. Il est facile de voir que c'est une sous-cat\'egorie tannakienne, stable par extension.  

\medskip $\bullet$ On doit \`a Christol et Mebkhout \cite{[CM3]},\cite[2]{[CM4]} la construction d'une filtration de
type Hasse-Arf sur $MCS(\sR/K)$: la filtration par les {\it pentes $p$-adiques}.

Le facteur $gr^{ \lambda}(M)$ est caract\'eris\'e par la propri\'et\'e suivante: pour tout $\rho <1$ suffisamment
proche de $1$, le rayon de convergence de toute solution de $gr^{ \lambda}(M)$ non nulle au point g\'en\'erique
$t_\rho$ de module $\rho$ est $\rho^{1+\lambda}$. 

Ou, ce qui revient au m\^eme: $F^{>\lambda}(M)=0 \iff $ $M$ admet une base de solutions analytiques dans le disque
non circonf\'erenci\'e de centre $t_\rho$ et de rayon $\rho^{1+\lambda}$ (pour tout $\rho $ suffisamment
proche de $1$).  

La hauteur du polyg\^one de Newton de $M$ s'appelle l'{\it irr\'egularit\'e $p$-adique} de $M$. C'est toujours un entier
(du fait qu'elle s'interpr\`ete comme un indice, fait beaucoup plus profond que dans le cas formel).  

\medskip $\bullet$  Soit ${E}=\bar K$ une cl\^oture alg\'ebrique de $K$, et posons $\sR_{\bar K}= \sR\otimes_K \bar
K$. On voit imm\'ediatement que $K$ est alg\'ebriquement ferm\'e dans $\sR$, donc $\sR_{\bar K}$ est int\`egre.

  Soit $MC(\sR_{\bar K}/\bar K)$ la cat\'egorie des $\sR_{\bar K}[{d\over dz}]$-modules de pr\'esentation
finie sur $\sR_{\bar K}$. Tout objet de $MC(\sR_{\bar K}/\bar K)$ provient par extension des scalaires d'un objet
de $MC(\sR\otimes_K L/L)$ pour une extension finie convenable $L/K$. On d\'eduit de ce qui pr\'ec\`ede qu'il est
libre sur $\sR_{\bar K}$. La cat\'egorie
$MC(\sR_{\bar K}/\bar K)$ est tannakienne sur
$\bar K$.

 On a une sous-cat\'egorie tannakienne $MCS(\sR_{\bar K}/\bar K)$ form\'ee des modules diff\'erentiels
solubles au bord ext\'erieur. La filtration par les pentes
$p$-adiques n'est pas sensible aux extensions finies des scalaires, donc se transporte sans probl\`eme en une
filtration de Hasse-Arf de $MCS(\sR_{\bar K}/\bar K)$.

\medskip $\bullet$ Soient $h$ un entier $\geq 1$, et $\sigma$ un automorphisme continu de $\bar K$ relevant
 l'automorphisme $x\mapsto x^{p^h}$ de $\bar k$. Il y a un unique automorphisme
$\sigma$-lin\'eaire $\varphi_h$ de $\sR=\sR_{K,z}$ tel que $z\mapsto z^{p^h}$. On dit qu'un objet $M$ de
$MC(\sR_{\bar K}/\bar K)$ {\it admet une structure de Frobenius} d'ordre $h$ (\resp une structure de Frobenius) si
$\varphi_h^\ast M\cong M$ (\resp pour $h$ non pr\'ecis\'e). 
 
D'apr\`es un th\'eor\`eme de Christol-Mebkhout (\cf \cite[6]{[CM4]}), la sous-cat\'egorie
pleine $MCF(\sR_{\bar K}/\bar K)$ de      
$MC(\sR_{\bar K}/\bar K)$ form\'ee des objets admettant une structure de Frobenius est en fait une sous-cat\'egorie
tannakienne, stable par extension (ce th\'eor\`eme devient d'ailleurs facile si on se limite \`a la
sous-cat\'egorie de $MC(\sR_{\bar K}/\bar K)$ form\'ee des objets semisimples, \cf \cite[5.3. rem. a]{an}). 
  Par un argument bien connu d\^u \`a Dwork, $MCF(\sR_{\bar K}/\bar K)$ est en fait une sous-cat\'egorie de
$MCS(\sR_{\bar K}/\bar K)$. 
Insistons sur le fait que {\it seule l'exis\-tence d'une structure
de Frobenius} est prise en compte dans la d\'efinition de $MCF(\sR_{\bar K}/\bar K)$, et non la structure de
Frobenius elle-m\^eme, \ie le choix d'un isomorphisme $\varphi_h^\ast M\cong M$. 

 \subsection{} En revanche, le formalisme des pentes frobeniusiennes attach\'ees aux $F$-cristaux ne rentre pas dans
le formalisme des filtrations par les pentes \'etudi\'e dans cet article: d'une part, les pentes peuvent \^etre
n\'egatives, d'autre part les conditions
$FP_4)$ et $FP_5)$ ne sont pas satisfaites. En fait, ces ``pentes" rentrent dans le formalisme de Saavedra
des $\otimes$-filtrations\footnote{comme me l'a fait observer J. Sauloy, les pentes rencontr\'ees dans la th\'eorie
des mo\-du\-les aux $q$-diff\'erences s'\'ecartent de m\^eme du formalisme de cet article (et s'apparentent en fait
aux pentes frobeniusiennes).}.   
\end{sloppypar}

\section{Induction et induction tensorielle}\label{s4}

\medskip
\subsection{Rappels, \protect \cf \cite[10.3]{[K4]}}\label{rappels}  
Soit $H$ un groupe. Soit $\Gamma\subset \mathfrak{S}_n$ un sous-groupe du
groupe des permutations de $\{1,2,\ldots , n\}$. Le {\it produit en
volute} $\Gamma \wr H$ est le produit semi-direct de $\Gamma$ avec $H^n$, $\Gamma$
agissant par permutation des facteurs; explicitement, $\sigma\in \Gamma\subset \mathfrak{S}_n$ agit par
\[\sigma^{-1}(h_1,\ldots,h_n)\sigma = (h_{\sigma(1)},\ldots,h_{\sigma(n)}). \] 

Supposons que $H$ soit sous-groupe d'indice fini $n$ d'un groupe $G$. Alors 
le choix d'un ensemble ordonn\'e $\gamma_1,\ldots,\gamma_n$ de repr\'esentants des classes \`a droite modulo $H$
d\'efinit un homomorphisme injectif
\[G\to \mathfrak{S}_n \wr H:\;\;g\mapsto \sigma_g.(\gamma_{\sigma_g(1)}^{-1}g\gamma_1,\ldots ,
\gamma_{\sigma_g(n)}^{-1}g\gamma_n), \]
o\`u $\sigma_g$ est la permutation d\'efinie par la condition que les $\gamma_{\sigma_g(i)}^{-1}g\gamma_i \in H$. 

Cet homomorphisme ne d\'epend qu'\`a automorphisme int\'erieur pr\`es du choix des $\gamma_i$.

\medskip Fixons un anneau de base $E$. Par $E$-repr\'e\-sen\-ta\-tion d'un groupe, nous sous-entendons que le
$E$-module sous-jacent est libre de type fini. Soit $V$ une $E$-repr\'e\-sen\-ta\-tion de $H$. Alors on munit 
$V^{  n}$ d'une
repr\'e\-sen\-ta\-tion du produit en volute $\mathfrak{S}_n \wr H$ en posant 
\[(\sigma^{-1}.(h_1,\ldots, h_n)) (v_1, \ldots , v_n)=(h_{\sigma(1)}v_{\sigma(1)}, \ldots ,
h_{\sigma(n)}v_{\sigma(n)} ).\]
Sa restriction \`a $G$ est une $E$-repr\'e\-sen\-ta\-tion de $G$, appel\'ee {\it induite} de $V$, et not\'ee
$Ind_H^G(V)$.

De m\^eme, on munit $V^{\otimes n}$ d'une repr\'e\-sen\-ta\-tion de $\mathfrak{S}_n \wr H$ en posant 
\[(\sigma^{-1}.(h_1,\ldots, h_n)) (v_1\otimes \ldots \otimes v_n)=h_{\sigma(1)}v_{\sigma(1)}\otimes \ldots \otimes
h_{\sigma(n)}v_{\sigma(n)} .\]
Sa restriction \`a $G$ est une $E$-repr\'e\-sen\-ta\-tion de $G$, appel\'ee {\it tenseur-induite} de $V$, et not\'ee
${}^{\scriptstyle\otimes}Ind_H^G(V)$.

Ces constructions sont compatibles \`a tout changement de base $E\to E'$, et sont fonctorielles en $V$.
 
\subsection{Quelques propri\'et\'es, \cf \cite[3.3, 3.15]{be}, \cite[10.3]{[K4]}} { $\;$}

\medskip 1) $\displaystyle Res_G^H[Ind_H^G(V)]\cong \oplus_{i=1}^{i=n}\;
{}^{\gamma_i}V,\;\;Res_G^H[{}^{\scriptstyle\otimes}Ind_H^G(V)]\cong \otimes_{i=1}^{i=n}\; {}^{\gamma_i}V$, 

o\`u ${}^{\gamma_i}V= \gamma_i\otimes V$ est la repr\'e\-sen\-ta\-tion conjugu\'ee de $V$ par $\gamma_i$.

\medskip 2) Si $H$ est normal dans $G$, on a un isomorphisme canonique \[Ind_H^G(V)\otimes_E Ind_H^G(W)\cong
\oplus_i\,Ind_H^G({}^{\gamma_i}V\otimes_E W)\]
  Cette propri\'et\'e se prouve \`a l'aide du th\'eor\`eme de d\'ecomposition de Mackey, \cf
\cite[3.3.5]{be}. En it\'erant, on trouve 

\[\displaystyle (Ind_H^G(V))^{\otimes n}\cong \bigoplus_{i_1,\ldots, i_{n-1}} \,Ind_H^G((\bigotimes_{j=1,\ldots n-1}
{}^{\gamma_{i_j}}V )\otimes V) \]
dont $Ind_H^G(Res_G^H({}^{\scriptstyle\otimes}Ind_H^G(V)))$ est un facteur direct. Si $n$ est inversible dans $E$,
${}^{\scriptstyle\otimes}Ind_H^G(V)$ est facteur direct de ce dernier $G$-module (\cite[3.6.9]{be}), donc facteur
direct de  $(Ind_H^G(V))^{\otimes n}$. 

\medskip 3) Soit $H'$ un sous-groupe d'indice $n$ d'un groupe $G'$. Soient $\phi: G\to G'$ un homomorphisme de
groupes induisant un homomorphisme $\psi: H\to H'$ et un {\it isomorphisme} $G/H \to G'/H'$. Choisissons des
repr\'esentants $\gamma_i, \gamma'_i$ compatibles qui se correspondent sous $\phi$. Soit $V'$ une
$E$-repr\'esen\-ta\-tion de $H'$, et notons
$\psi^\ast(V')$ la repr\'e\-sen\-ta\-tion de $H$ (de m\^eme module sous-jacent) via $\psi$.  

Alors on a des isomorphismes canoniques
\[Ind_H^G( \psi^\ast(V'))\cong \phi^\ast(Ind_{H'}^{G'}(
 V')),\;\;{}^{\scriptstyle\otimes}Ind_H^G( \psi^\ast(V'))\cong \phi^\ast({}^{\scriptstyle\otimes}Ind_{H'}^{G'}(
 V')).\]   
   
Cette propri\'et\'e se prouve en deux temps. Si $\phi$ est surjectif, le r\'esultat est imm\'ediat (voir aussi
\loccit 10.3.2.(4)). Si $\phi$ est injectif, cela r\'esulte du th\'eor\`eme de Mackey,
dans le cas d'une seule double classe (et de sa version multiplicative, \cf \loccit 10.3.3). Le cas g\'en\'eral
s'obtient par composition.

\subsection{Variante ``g\'eom\'etrico-alg\'ebrique"}\label{alg} La propri\'et\'e 3) ci-dessus jointe \`a la
compatibilit\'e de l'induction et de la tenseur-induction \`a tout changement de base $E\to E'$ permet de transposer
ces constructions du cadre des groupes abstraits au cadre des sch\'emas en groupes affines.  

Plus pr\'ecis\'ement, soient $G$ un sch\'ema en groupes affine sur $E$ et $H$ un sous-sch\'ema en
groupes ferm\'e tel que le quotient $G/H$ soit repr\'esentable par un $E$-sch\'ema fini constant de rang $n$ (\ie
par $(Spec \,E)^n$). Alors \`a toute $E$-repr\'e\-sen\-ta\-tion $V$ de $H$, on
peut attacher son induite $Ind_H^G(V)$ (\resp sa tenseur-induite
${}^{\scriptstyle\otimes}Ind_H^G(V)$) qui est une $E$-repr\'e\-sen\-ta\-tion de
$G$ de module sous-jacent $V^n$ (\resp $V^{\otimes n}$). Les propri\'et\'es ci-dessus valent encore dans ce cadre.

\medskip Dans toute la suite, ${E}$ sera un {\it corps de caract\'eristique nulle}.

\medskip \noindent Si $V$ est une re\-pr\'e\-sen\-ta\-tion semisimple de
$H$, $Ind_H^G(V)$ et ${}^{\scriptstyle\otimes}Ind_H^G(V)$ sont alors des re\-pr\'e\-sen\-ta\-tions
semisimples de $G$ (\cf \loccit 10.3.4).

 \subsection{Le cas normal} On suppose en outre que $H$ est normal dans $G$, et on identifie le groupe
quotient $G/H$ \`a un sous-groupe de $\mathfrak S_n$. On peut donc identifier
$G$
\`a un sous-groupe ferm\'e de $ G/H \wr  H\subset \mathfrak S_n\wr  H$. Notons que l'application $g\mapsto \sigma_g 
\in G/H \subset
\mathfrak{S}_n$ de \ref{rappels} n'est autre que l'action par translation \`a gauche de $G$ sur $G/H$. 

\medskip
On note aussi $\;\underline{H}\;$ l'image de $H\to GL(V)$,
$\;\bar G$ l'image de $G\to GL(V^{n})$ (re\-pr\'e\-sen\-ta\-tion induite), $\;\bar H$ l'image de $H$ dans $\bar G$, 
$\;{}^{\scriptstyle\otimes}\bar G$ l'image de $G\to GL(V^{\otimes n})$ (re\-pr\'e\-sen\-ta\-tion $\otimes$-induite), et enfin
$\;{}^{\scriptstyle\otimes}\bar H$ l'image de
$H$ dans ${}^{\scriptstyle\otimes}\bar G$. 

\begin{lemme}\label{L0} a) On a des plongements naturels \[\bar G\inj 
G/H \wr
\underline{H}  ,  \; {}^{\scriptstyle\otimes}\bar G\inj 
G/H \wr
\underline{H}.\] 
b) Si $\dim\,V\geq 1$ (\resp $\dim\,V>1$), ces plongements induisent des isomorphismes $\bar G/\bar H  
\cong G/H\,$ (\resp $\,\bar G/\bar H  
\cong {}^{\scriptstyle\otimes}\bar G/{}^{\scriptstyle\otimes}\bar H$);
\\ c) On a des triangles commutatifs naturels d'\'epimorphismes \[\begin{CD}
 G&&   \longrightarrow  && \bar G &  {\;\;\; } &  
 H&& \longrightarrow   && \bar H \\
& \displaystyle\searrow &&  \displaystyle\swarrow &{\;\;\;\;\;}&{\;\;\;\;\;}&
&   \displaystyle\searrow  && \displaystyle\swarrow &\\
 &&{}^{\scriptstyle\otimes}\bar G &&{\;\;\;\;\;}&
  {\;\;\;\;\;}&{\;\;\;\;\;} &{\;\;\;\;\;}& 
{}^{\scriptstyle\otimes}\bar H. 
\end{CD}\]
 \end{lemme} 
\prf a) est imm\'ediat. 

\noindent b) suit de ce que si $\dim\,V\geq 1$ (\resp $\dim\,V>1$), l'action de $\mathfrak{S}_n$ sur
$V^n$ (\resp $V^{\otimes n}$) par permutation des facteurs est fid\`ele. 

\noindent c) r\'esulte de ce que
${}^{\scriptstyle\otimes}Ind_H^G(V)$ est facteur direct de  $(Ind_H^G(V))^{\otimes n}$.\qed 

\subsection{Un cas particulier} Supposons $\underline{H}$ {\it simple}\footnote{\ie distinct de ${\mathbb G}_a$, et
n'admettant aucun sous-groupe ferm\'e normal distinct de $1$ et de lui-m\^eme. Un tel groupe alg\'ebrique est soit
fini, soit connexe semisimple.}, et supposons que $G/H$ soit le groupe cyclique $C_\ell$ \`a $\ell$ \'el\'ements,
avec
$\ell= n$ premier. Rappelons que dans le produit en volute
$C_\ell
\wr
\underline{H}$,
$C_\ell $ agit sur
$\underline{H}^\ell$ par permutation circulaire des facteurs
$\underline{H}$. Notons $p_i,
\,p_{ij}$ les projections de $\underline{H}^\ell$ sur le $i$-\`eme facteur (\resp sur le produit des $i$-\`eme et
$j$-\`eme facteurs). 

Commen\c cons par une variante du lemme de Goursat.

\begin{lemme}\label{L1} Soit $H'$ un sous-groupe ferm\'e de $C_\ell \wr \underline{H}$ contenu dans
$\underline{H}^\ell$, stable sous l'action int\'erieure de $C_\ell$, et tel que $p_1(H')=\underline{H}$. Alors
$H'=\underline{H}^\ell$ ou bien $H'$ est un groupe de la forme
\[\{(h=\phi_1(h),\phi_2(h),\ldots, \phi_{\ell}(h))\vert\,h\in \underline{H}\} \cong
\underline{H}\]  o\`u les $\phi_j$ sont des
automorphismes de $\underline{H}$. 
\end{lemme}

\prf La stabilit\'e de $H'$ sous $C_\ell $ entra\^{\i}ne imm\'ediatement que les
$p_j(H')=\underline{H}$ pour chaque $j$. D'apr\`es le lemme de Goursat, $p_{1j}$ est surjective \`a moins
que $p_{1j}(H')$ ne soit le graphe d'un automorphisme $\phi_j$. La stabilit\'e de $H'$ sous le groupe simple
$C_\ell $ entra\^{\i}ne que ceci arrive simultan\'ement pour tout $j\neq 1$, ou n'arrive pour aucun $j\neq 1$.
\qed

On suppose $V\neq 0$.

\begin{prop}\label{P1} a) Sous les hypoth\`eses pr\'ec\'edentes, on est dans l'un des deux cas suivants:
\\ Cas I) $\;\bar G=  C_\ell \wr \underline{H}$ (produit en volute); ses
seuls sous-groupes ferm\'es normaux sont 
 $\; 1,  \underline{H}^\ell,\bar G$. 
\\ Cas II)  $\;\bar G =C_\ell \ltimes \underline{H} $ (produit semi-direct); si le produit n'est pas isomorphe \`a un
produit direct, ses seuls sous-groupes ferm\'es normaux sont $\; 1,  \underline{H},\bar G$. 

Dans les deux cas, il n'y a pas de cha\^{\i}ne de sous-groupes ferm\'es normaux de $\bar G$ de longueur $>3$ ($1$ et
$\bar G$ \'etant compris).

 b) Si $\dim\,V >1$, on a $\;\bar G= {}^\otimes\bar G$. Dans le cas I), si $V$ est une repr\'e\-sen\-ta\-tion
absolument irr\'eductible de $H$, alors  ${}^{\scriptstyle\otimes}Ind_H^G(V)$ est une re\-pr\'e\-sen\-ta\-tion
 absolument irr\'educ\-ti\-ble de $G$ (et m\^eme de $H$). 
 \end{prop} 

\prf a) On sait que $\bar G/\bar H \cong C_\ell $ (lemme
\ref{L0}). Donc $\bar H$ est un sous-groupe ferm\'e de $\,C_\ell \wr \underline{H}\,$ stable par $C_\ell $. De
la formule $ Res_G^H[Ind_H^G(V)]\cong \oplus\; {}^{\gamma_i}V\;$, on d\'eduit que $\,p_1(\bar H)=\underline{H}$. On
d\'eduit du lemme pr\'ec\'edent qu'on est dans l'un des cas I) ou II). Soit $G'$ un sous-groupe ferm\'e normal non
trivial de $\bar G$. Alors $H'= G'\cap \underline{H}^\ell$ est un sous-groupe ferm\'e normal non trivial de $\bar H$
stable sous l'action int\'erieure de $C_\ell$. Comme $\underline{H}$ est simple, $p_1(H')=\underline{H}$ ou $1$. On
termine la classification des sous-groupes ferm\'es normaux en appliquant de nouveau le lemme pr\'ec\'edent.

 b) Montrons que l'on ne peut avoir simultan\'ement $\bar H=  
\underline{H}^\ell$ et ${}^\otimes\bar H= \underline{H} $, ou vice-versa. La seconde possibilit\'e est exclue
par le point $c)$ du lemme \ref{L0}: on a un \'epimorphisme $\bar H\surj {}^\otimes\bar H$. Ce point et le lemme
pr\'ec\'edent excluent en fait aussi la premi\`ere possibilit\'e: en effet, le noyau $\underline{H}^\ell \to
{}^\otimes\bar H=  \underline{H}$ serait un sous-groupe ferm\'e normal de $\underline{H}^\ell$ ($\neq
1,\neq\underline{H}^\ell$), stable par permutation circulaire des
$\ell$ facteurs, ce qui est impossible.
 
Pour la derni\`ere assertion de b), on remarque que si $V$ est une re\-pr\'e\-sen\-ta\-tion absolument irr\'eductible
de $H$ de dimension $>1$,  $\;\;Res_G^H[{}^{\scriptstyle\otimes}Ind_H^G(V)]\cong
\otimes\; {}^{\gamma_i}V$ est une re\-pr\'e\-sen\-ta\-tion absolument irr\'eductible de $ \underline{H}^\ell$.
\qed

Remarquons que comme $\underline{H}$ est simple, le groupe des automorphismes ext\'erieurs $Out(\underline{H})$ est
fini (si $\underline{H}$ est infini, cela d\'ecoule du th\'eor\`eme classique correspondant pour les alg\`ebres de
Lie semisimples).

\begin{lemme}\label{L2} Pla\c cons-nous dans le cas II) de la proposition pr\'ec\'edente. Supposons en outre que
$\ell$ ne divise pas $\vert Out(\underline{H})\vert $. Alors $V\cong
Res_G^H(W)= Res_{\bar G}^{\bar H}(W)$ pour une re\-pr\'e\-sen\-ta\-tion absolument irr\'eductible convenable $\,W$ de
$\bar G$, et
$\bar G$ est isomorphe au produit direct
$\, C_\ell
\times \underline{H} $.  
\end{lemme}

\prf L'hypoth\`ese entra\^{\i}ne que l'action $C_\ell \to
Aut(\underline{H})$ d\'efinissant le produit semi-direct $C_\ell \ltimes \underline{H} $ s'effectue par
automorphismes int\'erieurs, donc pro\-vient d'un homomorphisme $C_\ell \to \underline{H}$ puisque $\underline{H}$
est simple. On peut ainsi enrichir la re\-pr\'e\-sen\-ta\-tion $V$ en une re\-pr\'e\-sen\-ta\-tion $W$ de $\bar G$
(absolument irr\'eductible tout comme $V$) via
$C_\ell \to \underline{H}$, dont l'image n'est autre que $\underline{H}$. Quitte \`a composer cet homomorphisme
$\bar G \surj \bar H=\underline{H}$ avec un automorphisme de $\bar{H}$, on obtient une r\'etraction de $\bar H\inj
\bar G$.  
 \qed

 \section{Un th\'eor\`eme de structure pour les filtrations de Hasse-Arf}\label{s5}

Dans l'exemple classique des re\-pr\'e\-sen\-ta\-tions galoisiennes, l'efficacit\'e du formalisme des
groupes de ramification vient en grande partie de la possibilit\'e de faire une extension finie du corps de base
(lemme de Herbrand, etc\ldots). En particulier, pour une extension mod\'er\'ement ramifi\'ee $L/K$ d'indice de
rami\-fication $n$, on a $G_L^{(n\lambda)}=G_K^{(\lambda)}\cap G_L$.

Nous allons adapter une partie de ce formalisme de changement de base au cas des filtrations de type Hasse-Arf
abstraites.

\subsection{Cadre}\label{cadre}  
{$\;$}

$\bullet $ On fixe un ensemble ${\boP}$ de nombres premiers (\'eventuellement vide). On suppose que l'ensemble
${\boP}'$
 des nombres premiers
n'appartenant pas
\`a ${\boP}$ est {\it infini}. On pose $\hat{\Z}_{({\boP}')}=
\Pi_{\ell\notin {\boP}}\,\Z_\ell$.

\medskip $\bullet $  On se donne un groupe profini $\sI$. On suppose que $\sI$ s'inscrit dans une suite exacte 
\[\displaystyle{1\to \sP \to \sI \to \hat{\Z}_{({\boP}')}\to 1 }\; \] 
o\`u $\sP$ est un pro-${\boP}$-groupe \footnote{Une telle suite exacte se scinde pour raison de cardinal (\cf
\cite[5.9. cor. 1]{seCG}), mais peu importe ici.}. (En pratique,
$\sI$ sera le groupe d'inertie d'un corps henselien
$K$ pour une valuation discr\`ete \`a corps r\'esiduel
$k$ parfait, de sorte que $\vert \boP \vert = 0$ ou $1$; si $\boP$ est vide, $\sP$ est
trivial.) 
 
 Tout sous-groupe d'indice fini de $\hat{\Z}_{({\boP}')}$, en particulier l'image de tout sous-groupe ouvert $U$ de $\sI$,
est de la forme $n_U.\hat{\Z}_{({\boP}')}$, donc canoniquement isomorphe \`a $\hat{\Z}_{({\boP}')}$. Pour tout $n$
premier \`a ${\boP}$, la pr\'eimage de
$n\hat{\Z}_{({\boP}')}$ dans $U$ est not\'ee $U^{(n)}$, de sorte que $U/U^{(n)}\cong \Z/n\Z$. Noter que si $U$ est
normal dans $\sI$,
$U^{(n)}$ l'est aussi (comme intersection de $U$ et de la pr\'eimage de $nn_U\hat{\Z}_{({\boP}')}$ dans $\sI$).   

\medskip $\bullet $ Soit ${E}=\bar {E}$ un corps alg\'ebriquement clos de caract\'eristique nulle. On consid\`ere
$\sI$ comme ${E}$-sch\'ema en groupes affine constant, de m\^eme que ses sous-groupes ouverts.

 On se donne un \'epimorphisme de ${E}$-groupes affines
\[G\surj \sI.\] Pour tout sous-groupe ouvert normal $U\triangleleft \sI$, on pose $G_U= G\times_\sI U$. Pour toute
inclusion $U'\subset U$ de sous-groupes ouverts normaux de $\sI$, on note $\iota_{U'U}:
G_{U'}\subset G_U$ l'inclusion correspondante.  On a donc un syst\`eme inductif de cat\'egories tanna\-kiennes
neutralis\'ees
\[Rep_{E} G\to \ldots \to Rep_{E} G_U \stackrel{\iota_{{U'}U}^\ast= Res_{G_U}^{G_{U'}}}\to Rep_{E} G_{U'} \ldots  \] 
index\'e par les sous-groupes ouverts normaux de
$\sI$. On a par ailleurs des foncteurs $Ind_{G_{U'}}^{G_U}$ et ${}^\otimes Ind_{G_{U'}}^{G_U}$ qui vont dans l'autre
sens.

On peut identifier $Rep_{E} U$ (\resp $Rep_{E} \hat{\Z}_{({\boP}')}$) \`a la sous-cat\'egorie tanna\-kienne de
$Rep_{E} G_U$ form\'ee des re\-pr\'e\-sen\-ta\-tions $V$ telles que $\iota_{U'U}^\ast(V)$  (\resp
$\iota_{U^{(n)}U}^\ast(V)$) soit isomorphe \`a $\un^{\dim \,V}$ pour $U'$ convenable (\resp pour $n$ premier \`a
$\boP$ convenable). 

 \medskip $\bullet $ Enfin on suppose les $Rep_{E} G_U$ munies de {\it filtrations de Hasse-Arf} $(F^{>\lambda})$
avec les conditions de compatibilit\'e suivantes: pour tout $U$ (normal), tout $\lambda$, et tout $n$
premier
\`a ${\boP}$,

\medskip
  $Comp_1)\;\;$ $ F^{>n\lambda}\circ \iota_{U^{(n)},U}^\ast =\iota_{U^{(n)},U}^\ast\circ F^{>\lambda},$

\medskip
  $Comp_2)\;\;$ il existe une \'equivalence de cat\'egories tannakiennes \[\rho^\ast_{U^{(n)},U}:
Rep_{E} G_U
\cong Rep_{E} G_{U^{(n)}}\] compatible aux filtrations de Hasse-Arf\footnote{on n'impose aucune relation entre
$\rho^\ast_{U^{(n)},U}$ et 
$\iota^\ast_{U^{(n)},U}$.}. 

\begin{rem} En termes des filtrations $(G_U^{(\lambda)})$ attach\'ees aux filtrations de Hasse-Arf et par le
dictionnaire tannakien, ces conditions s'\'ecrivent de mani\`ere \'equivalente:

\medskip
$Comp_1')\;\;$ $\iota_{U^{(n)},U}(G_{U^{(n)}}^{(n\lambda)})= G_U^{(\lambda)}$, 

\medskip
  $Comp_2')\;\;$ il existe un isomorphisme $\rho_{U^{(n)},U}:\, G_{U^{(n)}} \cong G_U\;\;$ envoyant
$\,G_{U^{(n)}}^{(\lambda)}\,$ sur $\,G_U^{(\lambda)}$.

\medskip\noindent (Heuristiquement, ces deux conditions sont les
m\^achoires d'un
\'etau qui force l'existence de pentes fractionnaires, d\`es lors qu'il y a des pentes non nulles.) 
\end{rem}

\subsection{Le th\'eor\`eme de structure} Rappelons qu'un caract\`ere d'un ${E}$-groupe affine est un homomorphisme
de ce groupe affine vers $\bG_m$ (ou ce qui revient au m\^eme, une re\-pr\'e\-sen\-ta\-tion de dimension $1$ sur
${E}$).  Rappelons aussi que
le radical est le plus grand sous-groupe ferm\'e normal connexe pro-r\'esoluble.

Pour une repr\'esentation de $G$ ou $G_U$, nous emploierons de mani\`ere interchangeable les expressions ``{\it
mod\'er\'ee}" et ``{\it de pente $0$}". 

 \begin{thm}\label{struc} a) On fait les hypoth\`eses suppl\'ementaires suivantes sur les caract\`eres de $G_U$, pour
tout sous-groupe ouvert normal $U$ de $\sI$:

\medskip\noindent $\chi_1)$ tout caract\`ere d'ordre fini de $G_U$ provient d'un caract\`ere de $U$,
\\ $\chi_2)$ tout caract\`ere mod\'er\'e d'ordre fini de $G_U$ provient d'un caract\`ere du quotient
$\hat{\Z}_{({\boP}')}$ de $U$, et r\'eciproquement,
\\$\chi_3)$ toute re\-pr\'e\-sen\-ta\-tion irr\'eductible mod\'er\'ee de $G_U$ est de dimension un, \ie est un
caract\`ere. 

\medskip Alors les homomorphismes canoniques 
\[G/Rad(G) \to \pi_0(G)\to \sI\]
sont des isomorphismes. De plus le sous-groupe sauvage $G^{>0}$ est
extension du pro-${\boP}$-groupe $\sP$ par un groupe
pro-r\'esoluble connexe, et le quotient mod\'er\'e $G_{ 0}=G/G^{>0}$ est
extension du pro-${\boP}'$-groupe $\hat{\Z}_{({\boP}')}$ par un groupe
pro-r\'esoluble connexe. 
 
\medskip b) On suppose en outre que 

\medskip\noindent 
$Ext)\;\;$ pour tout $m$, il existe au plus une classe d'isomorphisme d'extensions it\'er\'ees
$N_m$ de $\un$ de dimension $m$ dans $Rep_{E} G$,
\\ $Ind)\;\;$ tout objet ind\'ecomposable de $Rep_{E} G$ est de la forme $V\otimes N_m$, o\`u $V$ est
irr\'eductible.

\medskip  Alors la composante neutre de $\, G^{>0}\,$ est un pro-tore, et $G_0$ est un groupe dont la composante
neutre est un pro-tore, ou bien le produit de $\bG_a$ par un tel groupe. 

\medskip c) Enfin, si de plus    
 
\medskip \noindent $\chi_4)$ tout caract\`ere de $G_U$ est d'ordre fini, 

\medskip alors $\;G^{>0}\cong \sP, \;$ et $\;\;G\cong \sI\,$ ou $\,\sI\times \bG_a$.
 \end{thm}  

 \prf a) Nous proc\'edons en cinq \'etapes, dont les principales sont 2) et 4).  

1) Examinons d'abord l'\'epimorphisme $\pi_0(G)\surj \sI$. Pour montrer que c'est un isomorphisme, il suffit de
montrer que pour toute re\-pr\'e\-sen\-ta\-tion $R$ o\`u $G$ agit \`a travers un groupe fini, $G$ agit en fait \`a
travers $\sI$. 

Notons $G^R$ le groupe fini $Im(G\to GL(R))$. Le syst\`eme des sous-groupes $G^R_U$ de $G^R$ (index\'e par le
syst\`eme ordonn\'e filtrant des sous-groupes normaux ouverts $U$ de $\sI$) est stationnaire. Fixons $U$ pour
lequel $G^R_U$ est mini\-mum. Il s'agit de voir que $G^R_U=1$. 

Supposons le contraire, et soit $\underline{H}$ un quotient simple de $G^R_U$. Alors la filtration de
$\underline{H}$ par les images des $G_U^{>\lambda}$ n'a qu'un saut, ce qui entra\^ine que les irr\'eductibles
$\neq \un $ de la sous-cat\'egorie
$Rep_{E} \underline{H}$ de
$Rep_{E} G_U$ sont purement d'une pente $\lambda$ (qui peut \^etre nulle). Nous allons montrer que $\underline{H}$
est ab\'elien.

Le cas de pente $0$ r\'esulte de l'hypoth\`ese $\chi_3)$. On va donc se restreindre au cas $\lambda\neq 0$. Par la
vertu des filtrations de Hasse-Arf, le nombre $\lambda.\dim V$ est un entier $>0$.

2) Si $\underline{H}$ est non-ab\'elien, il existe une re\-pr\'e\-sen\-ta\-tion
irr\'e\-duc\-ti\-ble $V$ de dimension $>1$; elle est fid\`ele puisque $\underline{H}$ est simple.

Soit $\ell$ un nombre premier $\in {\boP}'\,$ {\it ne divisant pas $\,\vert Out(\underline{H})\vert\,$, ni $\,
\dim\, V\;$ ni $\,\lambda.
\dim\, V\;$} (les ``c\^ot\'es" du polyg\^one de Newton). Il en existe puisque $\boP'$ est infini. 

Gr\^ace \`a $Comp_2)$, commen\c cons par transf\'erer le
probl\`eme sur
$U^{(\ell)}$ au moyen de $\rho^\ast_{U^{(\ell)},U}$. On voit donc
$\underline{H}$ comme l'image de $H=G_{U^{(\ell)}}$ dans $GL(V)$. 

Consid\'erons la re\-pr\'e\-sen\-ta\-tion induite $Ind_H^{G_U}(V)$. On est dans la situation d'application de
\ref{P1} et \ref{L2}, dont on reprend les notations (avec $G_U$ au lieu de $G$). Dans l'un ou l'autre des cas I),
II), il n'y a pas  de cha\^{\i}ne de sous-groupes ferm\'es normaux de $\bar G$ de longueur $>3$ ($1$ et
$\bar G$ \'etant compris). Il suit que la filtration de $\bar G$ par les images des $G^{>\lambda}$ au plus deux sauts. On a
$\iota_{U^{(\ell)},U}^\ast(Ind_H^{G_U}(V))\cong \oplus {}^{\gamma_i}V$, d'o\`u l'on d\'eduit par $Comp_1)$ que
$\lambda/\ell$ est l'un d'entre eux. Quant aux re\-pr\'e\-sen\-ta\-tions de $\bar G$ se factorisant par $C_\ell$,
leur restriction \`a $H$ est triviale; donc, toujours par $Comp_1)$, l'autre saut est $0$. Autrement dit, dans la
cat\'egorie tannakienne $Rep_{E} \bar G$ engendr\'ee par $Ind_H^{G_U}(V)$, on ne rencontre que les pentes
$\lambda/\ell$ et $0$; on a $\bar G^{(\lambda/\ell)}=\bar H$. 

Examinons le cas I) de \ref{P1}: l'objet ${}^\otimes Ind_H^{G_U}(V)$ de $Rep_{E} \bar G$ est alors irr\'eductible de
dimension $>1$. Par $\chi_3)$, sa pente est donc non nulle. Elle vaut donc
$\lambda/\ell$. Comme $\ell$ ne divise ni $\dim\, V$ ni $\lambda. \dim\, V$, donc pas non plus $\lambda.
\dim\, ({}^\otimes Ind_H^{G_U}(V))=\lambda. (\dim\, V)^\ell$, c'est impossible par \ref{int}.

Examinons le cas II): l'hypoth\`ese que $\ell
\not\vert \,Out(\underline{H})\vert$ entra\^{\i}ne que 
$V$ est de la forme $Res_{\bar G}^{\bar H}(W)$, pour un objet $W$ irr\'eductible de $Rep_{E} \bar G$ de dimension
$>1$,
\cf \ref{L2}. Par
$\chi_3)$, la pente de $W$ est donc non nulle. Elle vaut donc $\lambda/\ell$. Or c'est
impossible par \ref{int} puisque $\ell \not\vert \,\lambda. \dim\, W=\,\lambda. \dim\, V$. 

\medskip
3) On tire de ce qui pr\'ec\`ede que $\underline{H}$ est en fait cyclique d'ordre premier. Il en d\'ecoule que la
re\-pr\'e\-sen\-ta\-tion irr\'eductible $V$ est de dimension $1$. L'hypoth\`ese $\chi_1)$ contredit alors le
choix initial du sous-groupe $U$. On conclut que l'homomorphisme $\pi_0(G)\to \sI$ est un isomorphisme.

\medskip
 4) Examinons \`a pr\'esent l'\'epimorphisme $G/Rad(G)\surj \pi_0(G)$. Pour montrer que c'est un isomorphisme,
 il suffit de montrer que pour toute re\-pr\'e\-sen\-ta\-tion $R$ o\`u $G$ agit \`a travers $G/Rad(G)$, $G$ agit en
fait \`a travers un groupe fini.   
 
 D'apr\`es l'\'etape pr\'ec\'edente, on sait qu'il existe un sous-groupe ouvert normal $U \triangleleft \sI$ tel
que l'image de
$G_U$ dans $GL(\iota_{U,\sI}^\ast R)$ soit connexe. Quitte \`a remplacer $\sI$ par $U$, on peut donc supposer
l'image de $G$ dans $GL(R)$ connexe. Il s'agit de voir qu'elle est triviale.

Supposons le contraire, et soit $\underline{H}$ un quotient simple de cette image. Puisque $\underline{H}$ est
simple, la filtration de
$\underline{H}$ n'a qu'un cran, ce qui entra\^ine que les irr\'eductibles $\neq \un $ de la sous-cat\'egorie
$Rep_{E} \underline{H}$ de $Rep_{E} G$ sont purement d'une pente
$\lambda$ (qui ne peut \^etre nulle en vertu de $\chi_3)$). Comme $\underline{H}$ est non-ab\'elien (puisque
$Rad(G)$ agit trivialement), il existe une re\-pr\'e\-sen\-ta\-tion irr\'eductible $V$ de dimension $>1$; elle est
fid\`ele puisque $\underline{H}$ est simple.
 La suite de la d\'emonstration s'effectue comme en 2), et m\`ene \`a une contradiction. 

\medskip
5) On conclut de ce qui pr\'ec\`ede que la suite 
\[1\to  Rad(G)\to G \to \sI \to 1\] est exacte, ce qui prouve la premi\`ere partie de l'assertion $a)$. 

 Passons \`a la seconde partie de $a)$. Le radical du quotient mod\'er\'e $G_0$ n'est autre que l'image
du radical de $G$. On en d\'eduit que $G_0/Rad(G_0)$ s'identifie au quotient de $\sI$ par l'image $\sI^{>0}$
de $G^{>0}$ dans $\sI$. 
 
\noindent Toute repr\'esentation $V$ de $G$ provenant d'une repr\'esentation du quotient $G_0/Rad(G_0)$ est
semisimple et d'image finie. Par $\chi_3)$, elle est somme directe de caract\`eres mod\'er\'es de
$G$.  Soit $V$ une repr\'esentation irr\'eductible de
$G_0/Rad(G_0)$. Par $\chi_3)$, cette repr\'esentation provient d'un caract\`ere mod\'er\'e (d'ordre fini) de
$G$, qu'on peut identifier \`a des caract\`eres (d'ordre fini) de $\sI/\sI^{>0}$ d'apr\`es ce qui pr\'ec\`ede. On
conclut alors par $\chi_2)$ que $\sI/\sI^{>0} = \hat\Z_{(\boP')}$, donc que $\sI^{>0}= \sP$. Ceci ach\`eve la preuve
de a).
 
\medskip b) L'hypoth\`ese $Ext)$ entra\^ine que les sommes finies d'extensions it\'er\'ees de $\un$ forment
une sous-cat\'egorie tannakienne de $Rep_{E} G$ isomorphe \`a $Vec_{E}$ o\`u \`a $Rep_{E} \bG_a$. L'hypoth\`ese
$Ind)$ implique que toute re\-pr\'e\-sen\-ta\-tion de $G$ est somme directe finie d'une re\-pr\'e\-sen\-ta\-tion du quotient de $G$
par son radi\-cal unipotent, tensoris\'ee avec une extension it\'er\'ee de $\un$, d'o\`u le r\'esultat par le
dictionnaire tannakien, en combinant avec $a)$. 

\medskip c) Compte tenu de ce qui pr\'ec\`ede, il suffit de montrer que sous la condition $\chi_4)$, le radical de
$G$ co\"{\i}ncide avec le radical unipotent. Or ceci se teste sur les quotients de dimension
finie de $G$. On se ram\`ene comme au pas 4) au cas o\`u ces quotients sont connexes. L'assertion est alors
claire.
\qed

 \subsection{Premiers exemples} 1) Dans l'exemple \ref{gal} des re\-pr\'e\-sen\-ta\-tions galoisiennes d'un corps
strictement henselien, $Comp_1)$ et $Comp_2)$ sont sa\-tis\-faites de m\^eme que toutes les conditions du
th\'eor\`eme, mais ce dernier n'apporte rien: $G= \sI$. 

\medskip
2) Dans l'exemple \ref{ell} des repr\'esentations $\ell$-adiques locales, le th\'eor\`eme s'applique pour
${}^FG$ mais il ne donne rien de plus que le th\'eor\`eme de mono\-dromie locale $\ell$-adique de Grothendieck,
d'ailleurs
\'equivalent
\`a la condition $\chi_4)$ dans cette situation.

\medskip
3) Dans l'exemple \ref{levelt} des modules diff\'erentiels sur ${E}((z))$, on retrouve une partie des r\'esultats de
Katz \cite{[K1]}. Voici comment. 

On prend pour $\sI$ le groupe de Galois absolu de ${E}((z))$ (qui co\"{\i}ncide avec le groupe d'inertie
puisque ${E}$ est suppos\'e alg\'ebriquement clos; on a $\boP=\emptyset$). On a donc $\sI\cong \hat \Z$, et tout ouvert
$U$ est de la forme $\sI^{(n)}\cong \hat \Z$. 

Soit $\omega$ le foncteur fibre canonique de Katz
 sur $MC({E}((z))/{E})$ (\cf
\ref{levelt}), et posons $G= Aut^\otimes\omega$. Comme ce foncteur est construit comme la fibre en $1$ de l'extension
cano\-nique, on a un carr\'e commutatif pour tout $n>0$
\[\begin{CD}
MC({E}((z^{1/n}))/{E}) @>{\omega\;\cong\;}>>  Rep_{E} G_{\sI^{(n)}}\\
 @A \otimes_{{E}((z))}{E}((z^{1/n})) AA @A \iota_{\sI^{(n)}\sI}^\ast AA\\
   MC({E}((z))/{E}) @>{\omega\;\cong\;}>>  Rep_{E} G.
\end{CD}\]
La sous-cat\'egorie de $MC({E}((z))/{E})$ form\'ee des objets qui deviennent triviaux par ramification
kummerienne $z\mapsto z^n$ d'ordre $n$ non pr\'ecis\'e s'envoie sur $Rep_{E} \sI$, et corres\-pond, dualement, \`a un
\'epimorphisme $G\to \sI$. Il est bien connu qu'une telle ramification $\iota_{\sI^{(n)}\sI}^\ast$ multiplie les
pentes par $n$ (\cite[2.2.11.2]{[K1]}). Le changement de variable $z\mapsto z^{1/n}$ induit un isomorphisme
$MC({E}((z))/{E})\cong MC({E}((z^{1/n}))/{E})$ de cat\'egories tannakiennes compatible aux filtrations par les pentes, qui
se transporte en une \'equivalence $\rho^\ast_{\sI^{(n)},\sI}$. Enfin, les conditions $\chi_1),\chi_2),\chi_3), Ext),
Ind)$ du th\'eo\-r\`e\-me sont satisfaites (mais {\it pas $\chi_4)$}).

 Comme $\boP$ est vide, le th\'eor\`eme montre que
$G^{>0}$ est un pro-tore (connexe). Dans \cite[2.6.5]{[K1]}, Katz identifie son groupe des caract\`eres au groupe
ab\'elien ${E}((z))/{E}[[z]]$.

\section{Le contexte $p$-adique}\label{s7} 

\subsection{} Nous nous pla\c cons maintenant dans le cadre \ref{CM} des modules diff\'eren\-tiels sur des couronnes
$p$-adiques. Nous allons commencer par montrer comment il se moule dans le cadre \ref{cadre}.

 Pour nous raccrocher \`a la litt\'erature sur ce sujet
\cite{[F]}\cite{[M1]}\cite{[M2]}\cite{[T1]}\cite{[T2]}\cite{[T3]}, nous supposerons d\'esormais que {\it $K$ est \`a
valuation discr\`ete ($p$-adique), de corps r\'esiduel $k$ parfait}. 

Il d\'ecoule de l'hypoth\`ese de discr\'etion que le
sous-anneau
$\sE^\dagger=\sE^\dagger_{K,z}$ de
$\sR$ des fonctions qui sont born\'ees vers le bord ext\'erieur est en fait un corps. C'est m\^eme un corps
henselien (muni de la norme sup), de corps r\'esiduel $k((z))$. 

  On prend $\sI= G_{\bar k((z))}$, le groupe d'inertie du corps local $ k((z))$; $\sP$ est le groupe
d'iner\-tie sauvage. On a donc $\boP=\{p\}$. On prend $E=\bar K$.

\subsection{Changements de base finis} Toute extension finie s\'eparable de $k((z))$ est de
la forme $k'((z'))$ et se rel\`eve en une extension finie non ramifi\'ee
$\sE^\dagger_{K',z'}$ de
${\sE^\dagger_{K,z}}$ (on
prendra garde \`a ce que l'\'equation liant $z'$ \`a $z$ est en g\'en\'eral fort complexe; on peut toutefois choisir 
  $z'$ de sorte que $z$ soit un polyn\^ome sans terme constant en $z'$). 

\noindent La d\'erivation $d/dz$ et le Frobenius de $\sE^\dagger_{K,z}$ s'\'etendent (de mani\`ere unique)
\`a $\sE^\dagger_{K',z'}$ \`a (\cite[2.6]{[T3]}). Ils s'\'etendent aussi \`a l'extension \'etale finie
$\sR_{K',z'}:= \sR_{K,z}\otimes_{\sE^\dagger_{K,z}}\sE^\dagger_{K',z'}$ de
$\sR_{K,z}$. 
 
\medskip $\bullet$ Tout sous-groupe ouvert normal $U\triangleleft \sI$ correspond \`a une extension galoisienne
$\bar k((z_U))/ \bar k((z))$, qui donne lieu \`a une extension \'etale finie
$\sR_{\bar K,z_U}$ de
$\sR_{\bar K,z}$, galoisienne de groupe $\sI/U$ (et qu'on peut supposer contenue dans une clot\^ure alg\'ebrique
fix\'ee du corps des fractions de $\sR_{\bar K,z}$). Notons $\tilde\sR$ l'extension \'etale de $\sR_{\bar K,z}$ r\'eunion de ces
extensions finies. Elle est munie d'une action de $\sI$, d'un Frobenius \'etendant $\varphi$, et d'une d\'erivation
\'etendant $d/dz$ (elle d\'efinit en particulier un ind-objet de $MCF(\sR_{\bar K}/\bar K)$).

\medskip $\bullet$ Notons $MC(\tilde\sR/\bar K)$ la cat\'egorie tannakienne des modules diff\'erentiels de pr\'esentation finie sur
$\tilde \sR$. 
 Lorsque  $U'\subset U$ d\'ecrivent les sous-groupes ouverts normaux de $\sI$, les op\'erations \[\otimes_{\sR_{\bar
K,z_U  }}\sR_{\bar K,z_{U'} },\;\;\;\otimes_{\sR_{\bar K,z_U}}\tilde\sR  
 \] d\'efinissent une tour de
$\otimes$-foncteurs exacts fid\`eles
\[ MC(\sR_{\bar K,z_U}/\bar K)\overset{j^\ast_{U'U} }{\to}  MC(\sR_{\bar K,z_{U'}}/\bar K)\overset{j^\ast_{U'}
}{\to}  MC(\tilde\sR /\bar K)  ,\] qui se restreint en une tour de
$\otimes$-foncteur exacts fid\`eles 
\[ MCF(\sR_{\bar K,z_U}/\bar K)\overset{j^\ast_{U'U} }{\to}  MCF(\sR_{\bar K,z_{U'}}/\bar K) \overset{j^\ast_{U' }
}{\to} MCF(\tilde\sR /\bar K),\;
\hbox{\cf
\cite[3.2]{[T1]}}. \]  Ici $MCF(\tilde\sR /\bar K)$ d\'esigne la sous-cat\'egorie pleine de $ MC(\tilde\sR /\bar K)
$ form\'ee dans l'image d'un $j^\ast_{ U'}$. Remarquons que tout objet $\tilde M$ de $MC(\tilde\sR /\bar K)$ provient
d'un objet $M_U$ de $MC(\sR_{\bar K,z_U}/\bar K)$ pour $U$ convenable, donc est libre en tant que $\tilde R$-module
(puisque $M_U$ est libre sur $\sR_{\bar K,z_U}$, \cf \ref{CM}).

\medskip $\bullet$ On a de m\^eme, par adjonction, des foncteurs (exacts fid\`eles) dans l'autre sens: 
\[ MC(\sR_{\bar K,z_U'}/\bar K)\overset{j_{U'U\ast} }{\to}  MC(\sR_{\bar K,z_{U}}/\bar K)  ,\; MCF(\sR_{\bar
K,z_U'}/\bar K)\overset{j_{U'U\ast} }{\to}  MCF(\sR_{\bar K,z_{U}}/\bar K) .\]  

  $\bullet$ La cat\'egorie tannakienne des modules diff\'erentiels sur le corps de fractions $L=Q(\tilde \sR)$
admet un foncteur fibre \`a valeurs dans $Vec_{\bar K}$, \cf appendice \ref{neu}. On le fixe. On en d\'eduit un
foncteur fibre
\[\,\tilde \omega :  MCF(\tilde\sR /\bar K)\to Vec_{\bar K}\;,\;\] 
 puis des foncteurs fibre compatibles 
\[\,  \omega_U = \tilde\omega \circ j^\ast_{ U}: \; MCF(\sR_{\bar K,z_{U}}/\bar K)\to Vec_{\bar K}.\] 
On abr\`ege  $\omega_\sI$ en $\omega$, et on note $G$ le ${\bar K}$-sch\'ema en groupes
affine $ Aut^\otimes
\omega$.  

\medskip $\bullet$ Consid\'erons la sous-cat\'egorie tannakienne $MCF(\sR_{\bar K,z }/\bar K)_\sI$ de $MCF(\sR_{\bar
K,z }/\bar K)$ form\'ee des objets $M$ que
$j^\ast_{\sI}$ trivialise, \ie tels que l'homomorphisme canonique \[(M\otimes_\sR \tilde\sR)^\nabla\otimes_{\bar
K} 
\tilde\sR \to j^\ast_{\sI}(M)\] soit un isomorphisme. On a alors un isomorphisme canonique \[\omega(M)=
\tilde\omega((M\otimes_\sR \tilde\sR)^\nabla\otimes_{\bar
K} 
\tilde\sR )=(M\otimes_\sR \tilde\sR)^\nabla  \]
d'o\`u une action de $\sI$ sur $\omega(M)$ via son action sur $\tilde\sR$. On en d\'eduit aussi un homomorphisme 
$\sI\to  Aut^\otimes \omega_{\mid MCF(\sR_{\bar K,z }/\bar K)_\sI}$ (o\`u $\sI$ est vu comme ${\bar K}$-sch\'ema en
groupes affine). C'est en fait un isomorphisme, car $\omega_{\mid MCF(\sR_{\bar K,z }/\bar K)_\sI}$ (isomorphe au
foncteur fibre
$ \;M\mapsto (M\otimes_\sR \tilde\sR)^\nabla$) admet un quasi-inverse donn\'e par
 $V\mapsto  (\tilde\sR \otimes_{\bar
K} V)^{\sI}.$
 
 \medskip $\bullet$ L'inclusion $MCF(\sR_{\bar K,z }/\bar K)_\sI\subset MCF(\sR_{\bar
K,z }/\bar K)$ correspond donc, par le dictionnaire tannakien, \`a un \'epimorphisme 
\[G \surj \sI.\]
   L'action du groupe fini $\sI/U$ sur $\sR_{\bar K,z_U}$ induit une action de $\sI/U$
sur
$MCF(\sR_{\bar K,z_U}/\bar K)$ par auto-$\otimes$-\'equivalences. L'argument des droites de Katz
\cite[1.4.5]{[K1]} s'applique dans notre contexte\footnote{pour la commodit\'e du lecteur, rappelons bri\`evement en
quoi consiste cet argument. Le seul point \'epineux est de montrer que le noyau
$G_U$ de la fl\`eche $G=Aut^\otimes
\omega \to \sI/U $ est contenu dans $Aut^\otimes \omega_U$. Pour cela, il suffit d'apr\`es Chevalley de d\'emontrer
que toute droite $L$ dans une repr\'esentation arbitraire de $G$, qui est stable par
$Aut^\otimes \omega_U$, est en fait stable par $G_U$. Katz montre que $Aut^\otimes \omega$ permute les droites
distinctes conjugu\'ees de $L$ sous le groupe fini $\sI/U$, ce qui entra\^{\i}ne bien que $G_U$
stabilise ces conjugu\'ees, donc $L$.}, et montre que la suite 
\[ 1\to Aut^\otimes \omega_U \to Aut^\otimes
\omega \to \sI/U \to 1\] est exacte, d'o\`u 
 $Aut^\otimes \omega_U= G_U$. On a donc un diagramme strictement commutatif de $\otimes$-foncteurs
\[\begin{CD} MCF(\sR_{\bar K,z_U}/\bar K)&@>{\omega_U \;\cong}>> & Rep_{\bar K} G_U \\
j^\ast_{U\sI}\uparrow & &&&\uparrow \iota^\ast_{U\sI}\\
MCF(\sR_{\bar K,z}/\bar K)&@>{\omega  \;\cong}>> & Rep_{\bar K} G   \end{CD} \]   
o\`u les foncteurs $\iota^\ast_{U\sI}$ d\'esignent comme pr\'ec\'edemment les foncteurs de restriction des
repr\'esentations.

  De m\^eme, les foncteurs $j_{U\sI\ast}$ correspondent \`a l'induction des repr\'esentations.

\subsection{Filtration de Hasse-Arf sur les $Rep_{\bar K} G_U$ et compatibilit\'es en $U$}    

Ceci permet de d\'efinir la filtration de Hasse-Arf de $Rep_{\bar K}\, G_U$ en transportant la
filtration par les pentes $p$-adiques de $MCF(\sR_{\bar K,z_U},\bar K)$.

Examinons \`a pr\'esent le lien entre ces cat\'egories tannakiennes $ Rep_{\bar K} \,G_U$ (munies de leurs
filtrations de Hasse-Arf) pour une extension kummerienne.  Le changement de variable $z_U\mapsto z_{U^{(n)}}$ induit
un isomorphisme $\sR_{\bar K,z_U}\cong \sR_{\bar K,z_{U^{(n)}}}$, d'o\`u un isomorphisme
\[MCF(\sR_{\bar K,z_U}/\bar K)\cong MCF(\sR_{\bar K,z_{U^{(n)}}}/\bar K)\] de cat\'egories
tannakiennes compatible aux filtrations par les pentes, qui se transporte en une \'equivalence
$\rho^\ast_{U^{(n)},U}:\;Rep_{\bar K}G_{U}\cong Rep_{\bar K}G_{U^{(n)}}$ compatible aux filtrations de Hasse-Arf
(utiliser le sorite
\ref{sor}), d'o\`u
$Comp_2)$.  

De m\^eme, il est connu que $j^\ast_{U^{(n)}U}$ multiplie les pentes par $n$, \cf \cite[6.3.5]{[CM3]}. En utilisant de
nouveau le sorite \ref{sor}, on en d\'eduit que $\iota^\ast_{U^{(n)}U}$ fait de m\^eme, d'o\`u $Comp_1)$.

\subsection{V\'erification des hypoth\`eses de \ref{struc}}\label{rob} 

Les conditions $\chi_1),\chi_2)$ et $\chi_4)$ rel\`event de la th\'eorie de Robba des \'equations
diff\'erentielles de rang un (sur $\sR_{\bar K,z_U}$). 

Rappelons
bri\`evement les \'etapes (avec $U=\sI$ pour all\'eger les notations). Comme
$\sR^\ast=(\sE^\dagger)^\ast$ (et de m\^eme apr\`es toute extension finie de $K$), on voit que tout objet $M$ de
dimension
$1$ de
$MC (\sR_{\bar K}/\bar K)$ est d\'efini sur $\sE^\dagger_{L,z} $ pour une extension finie convenable $L/K$. Par
approximation-troncation \cf
\cite[44.5]{[CC]}, il est m\^eme d\'efini sur $L(z)$. La condition de solubilit\'e au bord ext\'erieur de la
couronne devient la condition de Dwork usuelle (solubilit\'e dans le disque g\'en\'erique $D(t_1,1^-)$). Un analogue
$p$-adique des produits de Weierstrass (\cf \cite[5.3, 5.4, 10.11]{R2})\footnote{Robba demande de passer \`a un
corps alg\'ebriquement clos et maximalement complet, mais Matsuda \cite{[M1]} a montr\'e qu'on peut s'en dispenser,
du moins pour $p>2$.} permet de trouver une base telle que l'op\'erateur diff\'erentiel correspondant \`a
$M$ s'\'ecrive
\[\sL = {d\over dz}+ {a_1\over z}+\ldots {a_n\over z^n}\;,
\;\;\;\vert a_i\vert
\leq 1\;\; ;\] o\`u $n-1=\lambda$ est la pente $p$-adique. La condition de solubilit\'e implique $a_1\in \Z_p$, mais
l'hypoth\`ese plus forte d'existence d'une structure de Frobenius implique\footnote{en fait: \'equivaut \`a cette
condition, \cf \cite{[CC]}; mais peu importe ici.}
$a_1\in \Z_{(p)}=\Z_p\cap \Q$.  
 
En particulier, $\sL\sim {d\over dz}+ {a_1\over z}$ si la pente est nulle, \ie si
 le caract\`ere correspondant de $ G$ est mod\'er\'e. On voit donc que l'homomorphisme r\'esidu, donn\'e par $a_1$
mod. $\Z$, induit un isomorphisme 
\[ X(G_{0})\cong\Z_{(p)}/\Z  \] o\`u $X(G_{0})$ d\'esigne le groupe des caract\`eres du quotient
mod\'er\'e de $G$. 

\noindent L'inverse de cet isomorphisme est d'ailleurs donn\'e par l'homomorphisme   
$X(\hat{\Z}_{(p')})\to X(G_{0})$ induit par $G_{0}\surj   \hat{\Z}_{(p')}$ (dualit\'e de
Pontriaguine). D'o\`u $\chi_2)$. 

 On a \[\sL^{\otimes p^N} = {d\over dz}+ {p^Na_1\over z}+\ldots {p^Na_n\over z^n}\sim {d\over dz}+ {p^Na_1\over z}\]
pour $N>>0$, donc le caract\`ere de $G$ correspondant \`a $\sL$ est fini. Ceci donne $\chi_4)$.

La condition $\chi_1)$ revient \`a dire que $\sL$ se trivialise sur
une extension $\sR_{\bar K,z_U}$ de $\sR_{\bar K,z}$ comme ci-dessus. Cela d\'ecoule de l'existence d'une structure
de Frobenius sur $\sL$, qu'on peut supposer d\'efinie sur $\sE^\dagger$ (puisque $\sR^\ast=(\sE^\dagger)^\ast$) et de
pente frobeniusienne nulle, \cf \cite[4.11]{[Cr1]} (Tsuzuki
\cite{[T3]} en donne une version en toute dimension).
 
  La condition $\chi_3)$ est nettement plus substantielle: c'est le ``{\it th\'eor\`eme de monodromie locale
$p$-adique mod\'er\'ee}" de  
\cite{[CM2]}, en pr\'esence d'une structure de Frobenius (voir aussi \cite{Dwork}).

 Enfin, les conditions $Ext)$ et $Ind)$ sont d\'emontr\'ees dans \cite[6.0.18, 6.0.16, 6.0.17]{[CM4]}
comme cons\'equence du th\'eor\`eme de l'indice de Christol-Mebkhout.

Nous sommes \`a m\^eme d'appliquer \ref{struc}, ce qui nous donnera le th\'eor\`eme de monodromie locale $p$-adique. 
 
\section{Le th\'eor\`eme de monodromie locale $p$-adique}\label{s8}   

\subsection{} Le th\'eor\`eme suivant est le r\'esultat principal de cet article. Il donne la structure de 
$MCF(\sR_{\bar K}/\bar K)\sim Rep_{\bar K}\,G$.      

On rappelle que le corps $p$-adique $K$ est suppos\'e
\`a valuation discr\`ete et \`a corps r\'esiduel $k$ parfait, que $\sI$
est le groupe d'inertie de $Gal(k((z))^{sep}/k((z)))$, et que $\sP$ est son $p$-Sylow (inertie sauvage). 
 
\begin{thm}\label{main} On a des isomorphismes canoniques \[G\cong \sI\times  \bG_a,\;\;
 G^{>0}\cong \sP \] de $\bar K$-groupes affines.
\end{thm}

\begin{compl} {\it L'isomorphisme $ \,G^{>0}\cong \sP\,$ identifie la filtration 
 $\;(G^{(\lambda)})$
provenant de la filtration par les pentes $p$-adiques \`a la filtration de $\sP$ par les groupes de ramification
sup\'erieurs. }   
\end{compl} 

\prf Le th\'eor\`eme r\'esulte du th\'eor\`eme \ref{struc} pour $U=\sI$, dont toutes les hypoth\`eses sont
remplies ainsi que nous l'avons vu au \S pr\'ec\'edent, et compte tenu du fait qu'il y a des extensions non triviales
de $\un$ par
$\un$. 

Le compl\'ement r\'esulte de l\`a et du th\'eor\`eme de Matsuda-Tsuzuki
\cite{[M2]}\cite{[T1]} (voir aussi \cite{[Cr2]} qui en donne une \'el\'egante preuve g\'eom\'etrique). Rappelons
la d\'emarche de Tsuzuki, l\'eg\`erement transpos\'ee \`a notre contexte: on a deux homomorphismes \[K_0(Rep_{\bar
K}\sI)\to \Z\] donn\'es l'un par le conducteur de Swan des repr\'esentations galoisiennes, l'autre par
l'irr\'egularit\'e
$p$-adique des objets de $MCF(\sR_{\bar K}/\bar K)$ associ\'es. Par
\ref{hauteur}, il s'agit de montrer que ces deux homomorphismes co\"{\i}ncident. Tensorisant avec $\Q$ et appliquant
le th\'eor\`eme d'induction d'Artin (\cf \cite[9]{se}), on se ram\`ene au
cas des ca\-rac\-t\`eres. En dimension un, le probl\`eme se traite par un calcul direct, voir \cite{[M1]},
\cite{[T1]}.
\qed

\begin{cor} Supposons $\;k\;$ fini. Alors les classes d'isomorphisme d'objets de 
$MCF(\sR_{\bar
K}/\bar K)$ forment un ensemble d\'enombrable.  
\end{cor} 
En effet, le groupe profini $\sI$ est alors m\'etrisable: l'ensemble de ses sous-groupes ouverts est d\'enombrable,
\cf
\cite[1.3]{seCG}.
\qed

\begin{ex} Pour tout premier $p$, les op\'erateurs de Bessel $B_\nu$ normali\-s\'es \`a la Dwork forment une famille
\`a un param\`etre $\nu \in \Z_p$ d'op\'erateurs diff\'e\-ren\-tiels solubles dans le disque g\'en\'erique. Vus
comme op\'erateurs \`a coefficients dans $\sE^\dagger$, ils sont en fait tous isomorphes  (et munis
d'une structure de Frobenius), \cf
\cite[7.5]{an}. Dans \loccit, on {\it d\'etermine} en fait la repr\'esentation galoisienne associ\'ee (le cas
int\'eressant est $p=2$), ce qui s'av\`ere nettement plus difficile que d'en prouver l'existence.  
\end{ex} 

  Passons de la structure de la cat\'egorie \`a la structure des objets. Voici deux traductions du
th\'eor\`eme, o\`u l'on redescend au corps de base
$K$. 

\begin{cor}  Tout objet de $MC(\sR/  K)$ admettant une structure de Frobenius a
une base de solutions dans $\,\sR_{K',z'}[\log z]\,$, o\`u $\,\sR_{K',z'}\,$ est l'extension finie \'etale de 
$\sR=\sR_{K,z}$ attach\'ee \`a une extension finie s\'eparable convenable de $k((z))$.  
\end{cor} 

(L'application directe du th\'eor\`eme ne donne le r\'esultat qu'apr\`es une extension finie $L/K'$, mais la
descente des constantes est standard, voir par exemple \cite[4.1.2]{[K1]})\footnote{au moment de la
r\'evision (Oct.01) de cet article soumis, l'auteur a pris connaissance d'une nouvelle pr\'epublication de K.
Kedlaya, intitul\'ee `a
$p$-adic local monodromy theorem', qui donne une autre d\'emonstration de l'\'enonc\'e 7.1.5. bas\'ee sur des
id\'ees compl\`etement diff\'erentes des n\^otres; il y est fait un usage intensif des structures de
Frobenius.}.\qed  
 
Nous dirons qu'un objet de $MC(\sR/  K)$ est absolument ind\'ecomposable (\resp absolument irr\'eductible) s'il le
reste
 apr\`es toute extension finie $L/K$ des scalaires.

\begin{cor} Tout objet {\rm absolument ind\'ecomposable} $M$ de $MC(\sR/  K)$ admettant une structure de Frobenius
devient, apr\`es une \'eventuelle extension finie du corps de base $K$, isomorphe \`a un
objet de la forme
\[M\otimes U_m\,,\]  o\`u $M$ provient d'une
$K$-re\-pr\'e\-sen\-ta\-tion (finie) absolument irr\'eductible de
$\sI$, et  $U_m\;$ est l'objet de dimension $m$ repr\'esent\'e par $\;(z{d\over dz})^m$ (extension it\'er\'ee
$m$-i\`eme de $\un$ par lui-m\^eme et ind\'ecomposable). En particulier, $M$ provient d'un objet de
$MC(\sE^\dagger/  K)$ admettant une structure de Frobenius.
\qed
\end{cor} 

De l\`a, on tire ais\'ement:

\begin{cor} Tout objet $M^\dagger$ de $MC(\sE^\dagger/  K)$ admettant une structure de Frobenius et tel que
$M^\dagger\otimes_{\sE^\dagger}\sR$ soit {\rm absolument irr\'eductible} a une base de solutions dans
l'extension finie non ramifi\'ee de 
$\sE^\dagger=\sE^\dagger_{K,z}$ attach\'ee \`a une extension finie s\'eparable convenable de $k((z))$.  \qed
\end{cor} 

Suivant \cite{[CM4]}, disons qu'un objet $M$ de $ MCS(\sR/  K)\;$ est {\it compl\`etement
irr\'eductible} si $gr^0(\underline{End}(M))=  K.id\;\;$. 

 \begin{cor} Tout objet compl\`etement irr\'eductible de $ MC(\sR/K)\;$ admettant une structure de Frobenius est
absolument irr\'eductible et de dimension une puissance de $\,p$.  
\end{cor} 

L'\'enonc\'e sur la dimension avait \'et\'e conjectur\'e 
dans \cite[3.0.12]{[CM4]} (en fait pour tout objet compl\`etement
irr\'eductible de pente $>0$ de $ MCS(\sR /  K)\;$).

\medskip
\prf Passons \`a l'objet correspondant $M$ de $ MCF(\sR_{\bar K}/\bar K)\;$, et \`a la
repr\'esentation $V$ de $G$ associ\'ee. On a $\dim M=\dim_{\bar K} V$. L'hypoth\`ese de compl\`ete
irr\'eductibilit\'e
 se traduit par \[(End(V))^{G^{>0}}=\bar K.id,\] c'est-\`a-dire par
l'irr\'eductibilit\'e de
$\underline{End}(Res_{G}^{G^{>0}}(V))$ d'apr\`es Schur (noter que puisque $G^{>0}$ agit \`a travers un
groupe fini, cette repr\'esentation est semisimple). A fortiori,
$Res_{G}^{G^{>0}}(V)$ est irr\'eductible (donc $V$ aussi). Comme en vertu du th\'eor\`eme \ref{main},
$G^{>0}=\sP$ est un pro-$p$-groupe, la dimension de cette repr\'esentation est une puissance de $p\,$
(\cite[6.5]{se}).\qed

\subsection{Sur une conjecture de Dwork} L'existence de structures de Frobenius n'est que faiblement utilis\'ee dans
le texte. Il serait int\'eressant de remplacer $MCF(\sR_{\bar K}/\bar K)$ par la sous-cat\'egorie tannakienne
maximale
$MCS^{\Z_{(p)}/\Z}(\sR_{\bar K}/\bar K)$ de 
 $ MCS(\sR_{\bar K}/\bar K)$ form\'ee des objets d'exposants et
r\'esidus d\'eterminantiels non-Liouville et dans $\;{\Z_{(p)}/\Z}$, \cf \cite[7]{[CM4]}. On obtiendrait ainsi la
co\"{\i}ncidence de $\;MCF(\sR_{\bar K}/\bar K)\;$ et 
$\;MCS^{\Z_{(p)}/\Z}(\sR_{\bar K}/\bar K)$ conjectur\'ee dans
\loccit (et ant\'erieurement, sous une forme l\'eg\`erement moins pr\'ecise, par Dwork, comme me l'a fait remarquer
G. Christol).  Le seul probl\`eme qu'on rencontre pour ce faire est celui de la stabilit\'e de
$MCS^{\Z_{(p)}/\Z}(\sR_{\bar K}/\bar K)$ lorsqu'on passe
\`a une extension finie \'etale de $\sR_{\bar K}$ induite par une extension finie s\'eparable de $\bar k((z))$ du
type Artin-Schreier. Certains modules diff\'erentiels peuvent acqu\'erir des pentes nulles lors d'une telle
extension, et il ne semble pas facile alors d'en contr\^oler les exposants.  

Pour tourner (quelque peu artificiellement) cette difficult\'e, on peut introduire la sous-cat\'egorie pleine
$MCS^{\Z_{(p)}/\Z}_{stable}(\sR_{\bar K}/\bar K)$ de
$MCS^{\Z_{(p)}/\Z}(\sR_{\bar K}/\bar K)$ form\'ee des objets $M$ tels que pour tout sous-groupe normal
ouvert $U\triangleleft
\sI$, $j_{U\sI}^\ast(M)$ appartienne \`a $MCS^{\Z_{(p)}/\Z}(\sR_{{\bar K},z_U}/\bar K)$. Comme $j_{U\sI}^\ast$ est
un $\otimes$-foncteur, $MCS^{\Z_{(p)}/\Z}_{stable}(\sR_{\bar K}/\bar K)$ est bien une sous-cat\'egorie tannakienne
de $MC(\sR_{\bar K}/\bar K)$. Elle contient $MCF(\sR_{\bar K}/\bar K)$. On dispose de foncteurs \[j_{U\sI}^\ast:
MCS^{\Z_{(p)}/\Z}_{stable}(\sR_{\bar K}/\bar K)\to MCS^{\Z_{(p)}/\Z}_{stable}(\sR_{{\bar K},z_U}/\bar K) . \]

On peut alors remplacer dans tout ce qui pr\'ec\`ede la cat\'egorie $\;MCF(\sR_{\bar K}/\bar K)\;$ par
$\;\;MCS^{\Z_{(p)}/\Z}_{stable}(\sR_{\bar K}/\bar K)$\footnote{pour prouver $\chi_1)$, on dispose du th\'eor\`eme
de Chiarellotto-Christol \cite{[CC]} qui montre que les objets de dimension $1$ de
$MCS^{\Z_{(p)}/\Z}_{stable}(\sR_{\bar K}/\bar K)$ et de $MCF(\sR_{\bar K}/\bar K)$ co\"{\i}ncident.}, et $\,G \,$
par le groupe tannakien correspondant. Au bout du compte, on obtient la variante affaiblie suivante de la conjecture
de Dwork:

\begin{scolie} $MCS^{\Z_{(p)}/\Z}_{stable}(\sR_{\bar K}/\bar
K)=MCF(\sR_{\bar K}/\bar K).$\qed
\end{scolie}

En termes heuristiques: pour les modules diff\'erentiels solubles au bord d'une mince couronne, imposer la
rationalit\'e des exposants
$p$-adiques sous une forme suffisamment forte entra\^{\i}ne l'existence d'une structure de Frobenius. 
  
\subsection{Extensions et foncteurs fibres canoniques} 

Katz \cite{[K0]} a introduit la notion d'extension canonique d'une extension s\'eparable de $k((z))$ \`a un
rev\^etement fini \'etale de $\mathbb{G}_{m, k}$. Cela fournit une \'equivalence de cat\'egories pourvu qu'on se
li\-mi\-te aux rev\^etements `sp\'eciaux' de $\mathbb{G}_{m, k}$, c'est-\`a-dire mod\'er\'es \`a l'infini, et dont le
groupe de monodromie g\'eom\'etrique a un unique $p$-Sylow.  
Cette correspondance a \'et\'e transport\'ee
au cadre des modules diff\'erentiels quasi-unipotents par Matsuda
\cite{[M2]}.   

 En fait, le corollaire \ref{main} permet de d\'efinir, plus directement, la notion d'extension canonique de tout
objet de $MCF(\sR_{\bar K}/\bar K)$ en un $\bar K$-isocristal surconvergent sur $\mathbb{G}_{m, \bar k}$ (muni d'une
structure de Frobenius): on peut supposer l'objet ind\'ecomposable, et avec les
notations de \loccit, l'extension de $U_m$ est \'evidente, et celle de $M$ s'obtient \`a partir de l'extension de
Katz de la repr\'esentation finie de l'inertie correspondante (\cf \cite{[M2]}, nous laissons les d\'etails au
lecteur).

\noindent En prenant la fibre de l'extension canonique au point $1\in \mathbb{G}_{m, \bar k}$, on obtient alors un
foncteur fibre canonique sur
$MCF(\sR_{\bar K}/\bar K)$.
 
\newpage 
\renewcommand{\appendixname}{Appendice}

\appendix

\section{Sur les dimensions de re\-pr\'e\-sen\-ta\-tions irr\'eductibles}\label{s9}

\medskip
\subsection{} On suppose encore que ${E}$ est un corps de caract\'eristique nulle.

\begin{prop}\label{Papp} Soit $G$ un groupe alg\'ebrique lin\'eaire sur ${E}$. On suppose que toute re\-pr\'e\-sen\-ta\-tion
irr\'eductible de $G$ l'est absolument. Alors les dimensions de ses re\-pr\'e\-sen\-ta\-tions
irr\'eductibles non triviales sont premi\`eres entre elles dans leur ensemble.  
\end{prop}

\prf L'hypoth\`ese entra\^{\i}ne que toute re\-pr\'e\-sen\-ta\-tion irr\'eductible de $G$ sur $\bar{E}$ est d\'efinie sur
${E}$. On peut donc remplacer
${E}$ par sa cl\^oture alg\'ebrique. On distingue plusieurs cas:

\noindent a) si $G$ est non connexe, le groupe (fini) des composantes connexes $\pi_0(G)$ admet des re\-pr\'e\-sen\-ta\-tions
irr\'eductibles non triviales. Quitte \`a remplacer $G$ par $\pi_0(G)$, on peut donc supposer $G$ fini. Soient $n_i$
les dimensions des repr\'e\-sen\-ta\-tions
irr\'eductibles de $G$ (avec $n_0= 1$ pour la repr\'esen\-ta\-tion triviale). On sait d'une part que $n_i$ divise
$\vert G\vert$ (\cite[6.5]{se}), et d'autre part que $\vert G\vert = \sum n_i^2 = 1+ \sum_{i>0} n_i^2$
(\cite[2.4]{se}). Ces deux faits entra\^{\i}nent imm\'ediatement que le pgcd des $n_i,\,i>0$ est \'egal \`a $1$
\footnote{signalons \`a ce propos un r\'esultat moins \'el\'ementaire de Thompson \cite{th}: si les dimensions
$n_i>1$ sont toutes divisibles par $p$, alors $G$ est $p$-nilpotent. De nombreux r\'esultats dans cette
direction ont \'et\'e obtenus depuis par l'\'ecole de Mayence.}. 

Supposons \`a pr\'esent $G$ connexe. 

\noindent b) Si $G$ est r\'esoluble, toute re\-pr\'e\-sen\-ta\-tion irr\'eductible est de dimension $1$.

\noindent c) Si $G$ n'est pas r\'esoluble, le groupe semisimple quotient de $G$ par son ra\-di\-cal admet des
re\-pr\'e\-sen\-ta\-tions irr\'eductibles non triviales. Quitte \`a remplacer $G$ par ce quotient, on peut donc supposer $G$
semisimple. On note comme d'habitude
$\rho$ la demi-somme des racines positives. Soit $V(2m\rho),
\,m>0$ la re\-pr\'e\-sen\-ta\-tion
irr\'eductible de $G$ de plus haut poids $2m\rho$ (a priori c'est une
re\-pr\'e\-sen\-ta\-tion de
$Lie\,G$, mais elle s'int\`egre en une re\-pr\'e\-sen\-ta\-tion de $G$ puisque $2m\rho$ est dans le r\'eseau des racines). 
La formule des caract\`eres de Weyl montre que la
dimension de $V(2m\rho)$ est
$(2m+1)^N$, o\`u $N$ est le nombre de racines positives,
\cf \cite[\S9, ex. 1]{lie}. Ces dimensions sont premi\`eres entre elles. \qed

\begin{cor}\label{corapp} Soit $\sT$ une cat\'egorie tannakienne sur ${E}$, munie
d'une filtration de Hasse-Arf. On suppose que $\sT$ contient des objets qui ne sont pas extension it\'er\'ee de
$\un$. On suppose qu'il n'y a aucune pente enti\`ere non nulle (ce qui est en particulier
le cas si les pentes sont
$<1$).
 Alors $\sT$ contient des objets irr\'eductibles $\not\cong \un$ de pente nulle.
\end{cor} 

\prf On peut supposer $\sT$ alg\'ebrique, \ie ayant un g\'en\'erateur tensoriel. Quitte \`a remplacer $F$ par une
extension finie (\cf
\ref{sous/ext}), on peut supposer $\sT$ neutralis\'ee, et que tout objet irr\'eductible l'est
absolument. On peut m\^eme supposer, en rempla\c cant le groupe tannakien associ\'e
$G$ par un quotient (et $\sT$ par la sous-cat\'egorie
tannakienne correspondante) que $G$ est simple. Supposons par l'absurde que tout objet de pente
nulle soit trivial, \ie somme de copies de $\un$. Alors $G$ n'a qu'une pente, qui est non enti\`ere par
hypoth\`ese. Toute re\-pr\'e\-sen\-ta\-tion irr\'eductible non triviale a une dimension divisible par le d\'enominateur de
cette pente, ce qui est absurde d'apr\`es la proposition pr\'ec\'edente.\qed

\begin{rem} Ceci s'applique \`a la cat\'egorie tannakienne engendr\'ee par le module diff\'erentiel de Bessel
sur l'anneau de Robba $\sR$ diadique, \cf \cite{an}. Dans ce cas, les pentes sont $0, 1/3, 1/2$. On trouve des
irr\'eductibles non-triviaux de pente nulle, mais trivialisables par ramification kummerienne $z\mapsto z^3$.
\end{rem}

\renewcommand{\appendixname}{Appendice}

 \section{Extensions de Picard-Vessiot infinies}\label{s10}

\subsection{} Dans cet appendice, on \'etend, faute de r\'ef\'erence, une petite partie de la th\'eorie de
Picard-Vessiot au cas d'une cat\'egorie tannakienne non n\'ecessairement alg\'ebrique de modules diff\'erentiels. 

Soit $(L,\partial )$ un corps diff\'erentiel de corps de constantes alg\'ebriquement clos $E=L^\partial$. Par
module diff\'erentiel sur $L$, on entend un $L[\partial]$-module de dimension finie en tant que $L$-espace. Si
$(\sL,\partial)$ est un corps diff\'erentiel extension de  $(L,\partial )$, on dit qu'un module diff\'erentiel $M$
sur $L$ est soluble dans $\sL$ si l'homomorphisme canonique 
$ (M\otimes_L \sL)^\partial\otimes_E \sL \to M\otimes_L \sL $
est un isomorphisme.
 
\medskip\noindent Soit $\sT$ une cat\'egorie tannakienne sur $E$ dont les objets sont des modules diff\'erentiels sur
$L$. Observons que $\sT$ est essentiellement petite, et que le cardinal de l'ensemble des classes d'isomorphisme
d'objets est $\leq \vert L\vert$.

\begin{defn} On dit que $(\sL, \partial)$ est une extension de
Picard-Vessiot (fractionnaire) de $(L,\partial)$ pour $\sT$ si les conditions suivantes sont satisfaites:  
  \\ $i)$ $\sL^\partial= E$,
\\ $ii)$ tout objet de $\sT$ est soluble dans $\sL$,
\\ $iii)$ $(\sL, \partial)$ est minimal pour ces propri\'et\'es.  
\end{defn}

Lorsque $\sT$ est alg\'ebrique (\ie si tout objet est sous-quotient d'une somme finie de $M^{\otimes n}\otimes
(M^\vee)^{\otimes m}$ pour un objet $M$ convenable), l'existence - et l'unicit\'e \`a isomorphisme non unique pr\`es
- d'une extension de Picard-Vessiot est connue, \cf  \eg \cite[3,
A.1]{vdp}, \cite[3.4.4.4]{an0}. 
 
\medskip Sous $ii)$, il est facile de voir que la condition $iii)$
\'equivaut \`a 

\medskip\noindent $iv)$ $\sL$ est engendr\'e, en tant que corps, par les sous-$L$-espaces\footnote{ces
sous-espaces sont aussi naturellement des modules diff\'erentiels.}
$<N^\vee, (N\otimes_L \sL)^\partial>$, o\`u $N$ d\'ecrit un ensemble de repr\'esentants des classes d'isomorphismes
d'objets de
$\sT$, et
$<\;>$ est le crochet de dualit\'e.
 
\begin{defn} Nous dirons que $(\sL, \partial)$ est une extension de
Picard-Vessiot partielle de $(L,\partial)$ pour $\sT$ si les conditions $i)$ et $iv)$ ci-dessus sont
satisfaites. 
\end{defn}

En vertu de $iv)$ et du fait que le cardinal de l'ensemble des classes d'isomorphisme
d'objets de $\sT$ est $\leq \vert L\vert$, toute extension de
Picard-Vessiot partielle v\'erifie $\vert \sL\vert =\vert L\vert$.
   
 \begin{prop} Il existe une extension de Picard-Vessiot $(PV(\sT), \partial)$ de $(L,\partial)$ pour $\;\sT$.
 \end{prop}

\prf  Fixons un ensemble $\Omega$ contenant
$L$ et de cardinal strictement sup\'erieur
\`a celui de $L$. Consid\'erons l'ensemble $\sP$ des
quadruplets $(L', +, . , \partial)$ o\`u $L'$ est un sous-ensemble de $\Omega$ contenant $L$, o\`u $+,.$ sont des
lois de composition faisant de $L'$ un corps, et o\`u $\partial$ est une d\'erivation faisant de $(L',\partial)$ une
extension de Picard-Vessiot partielle de $(L,\partial)$ pour $\sT$. On munit $\sP$ d'un ordre partiel
corres\-pondant \`a la notion de sous-corps diff\'erentiel. Il est facile de voir que la r\'eunion d'une cha\^ine
d'\'el\'ements de $\sP$ est encore un \'el\'ement de $\sP$. Par Zorn, il existe un \'el\'ement maximal contenant $
(L, +, . , \partial)$. Notons $(PV(\sT), \partial)$ ce corps diff\'erentiel.

Il s'agit de montrer que tout objet $M$ de $\sT$ est soluble dans $PV(\sT)$. Pour cela consid\'erons une extension
de Picard-Vessiot $\sL$ de $ PV(\sT)$ pour (la cat\'egorie
tannakienne engendr\'ee par) le module diff\'erentiel $M\otimes_L  PV(\sT)$ sur $PV(\sT)$. On a $\vert
\sL\vert=\vert L\vert$, donc $\sL$ est plongeable dans $\Omega$. 

Il suffit alors de montrer que
$\sL= PV(\sT)$, ce qui, par maximalit\'e de $PV(\sT)$, revient \`a montrer que $\sL$ est une extension de
Picard-Vessiot partielle de $L$ pour $\sT$. 
 
V\'erifions $i)$: comme $\sL$ est une extension de Picard-Vessiot pour $M\otimes_L  PV(\sT)$, on a $\sL^\partial=
PV(\sT)^\partial$, donc
$\sL=E$.

V\'erifions $iv)$: soit $\sL'$ le sous-corps de $\sL$ engendr\'e par les $<N^\vee, (N\otimes_L\sL)^\partial>$, o\`u
$N$ d\'ecrit un ensemble de repr\'esentants des classes d'isomorphismes d'objets de
$\sT$. Il contient $PV(\sT)$, donc aussi les $<N^\vee\otimes_L PV(\sT), (N\otimes_L\sL)^\partial>$  o\`u $N$
d\'ecrit  un ensemble de repr\'esentants des classes d'isomorphismes d'objets de la cat\'egorie
tannakienne engendr\'ee par le module diff\'erentiel $M\otimes_L  PV(\sT)$. Comme $\sL$ est une extension de
Picard-Vessiot pour $M\otimes_L  PV(\sT)$, on conclut que $\sL'=\sL$.
  \qed
  
\begin{cor}\label{neu} La cat\'egorie tannakienne $\sT$ est neutre.
\end{cor}

En effet, comme dans le cas usuel (o\`u $\sT$ est alg\'ebrique), on a un
foncteur fibre $\sT \to Vec_E$ ``donn\'e par"
$M\mapsto (M\otimes_L PV(\sT))^\partial$. \qed

\begin{rem} Comme dans le cas usuel, on peut montrer que $PV(\sT)$  est unique \`a isomorphisme non unique pr\`es.
\end{rem}

 \end{document}